\documentclass[a4paper]{amsart}
\usepackage{a4wide}
\usepackage{amssymb}
\usepackage[latin1]{inputenc}
\usepackage[all]{xy}

\title{Unimodular $L_\infty$-algebras}
\author{Johan Granåker}
\address{Department of Mathematics \\ Stockholm University \\ 106 91 Stockholm}
\email{johang@math.su.se}
\renewcommand{\k}{\mathrm{l\hspace{-1.5pt}k}}
\newcommand{\w}{\mathbf{w}}
\newcommand{\cC}{\mathcal{C}}
\newcommand{\cE}{\mathcal{E}}
\newcommand{\cF}{\mathcal{F}}

\newcommand{\cO}{\mathcal{O}}
\newcommand{\cP}{\mathcal{P}}
\newcommand{\cT}{\mathcal{T}}
\newcommand{\fG}{\mathfrak{G}}
\newcommand{\fg}{\mathfrak{g}}
\newcommand{\one}{\mathrm{1}\hspace{-2.5pt}\mathrm{l}}
\newcommand{\ra}{\longrightarrow}
\newcommand{\cal}{\circlearrowleft}

\newcommand{\isoto}{\xrightarrow{\sim}}
\DeclareMathOperator{\tr}{tr}
\DeclareMathOperator{\ad}{ad}

\DeclareMathOperator{\Lie}{\mathcal{L}\mathit{ie}}
\DeclareMathOperator{\Com}{\mathcal{C}\mathit{om}}
\DeclareMathOperator{\Ass}{\mathcal{A}\mathit{ss}}
\DeclareMathOperator{\ULie}{\mathcal{U}\!\mathcal{L}\mathit{ie}}
\DeclareMathOperator{\UAss}{\mathcal{U}\!\mathcal{A}\mathit{ss}}
\DeclareMathOperator{\sgn}{sgn}
\DeclareMathOperator{\Hom}{Hom}
\DeclareMathOperator{\End}{End}
\DeclareMathOperator{\Der}{Der}
\DeclareMathOperator{\Ber}{Ber}

\DeclareMathOperator{\Aut}{Aut}
\DeclareMathOperator{\Bij}{Bij}
\DeclareMathOperator{\oEnd}{\mathcal{E}\!\mathit{nd}}
\DeclareMathOperator{\id}{id}

\DeclareMathOperator{\cone}{cone}
\DeclareMathOperator{\Cyc}{Cyc}
\renewcommand{\div}{\mathrm{div}}
\newcommand{\iso}{\cong}
\newcommand{\decor}[2]{#1\langle#2\rangle}
\theoremstyle{plain}
\newtheorem{thm}{Theorem}
\newtheorem{proposition}{Proposition}
\newtheorem{lemma}{Lemma}
\newtheorem{corollary}{Corollary}
\theoremstyle{definition}
\newtheorem*{definition}{Definition}
\newtheorem*{remark}{Remark}
\newtheorem*{example}{Example}
\newcommand{\bde}{\begin{definition}}
\newcommand{\ede}{\end{definition}}
\newcommand{\bth}{\begin{thm}}
\renewcommand{\eth}{\end{thm}}
\newcommand{\bpr}{\begin{proposition}}
\newcommand{\epr}{\end{proposition}}
\newcommand{\ble}{\begin{lemma}}
\newcommand{\ele}{\end{lemma}}
\newcommand{\bpro}{\begin{proof}}
\newcommand{\epro}{\end{proof}}
\newcommand{\bre}{\begin{remark}}
\newcommand{\ere}{\end{remark}}
\newcommand{\bex}{\begin{example}}
\newcommand{\eex}{\end{example}}
\newcommand{\bco}{\begin{corollary}}
\newcommand{\eco}{\end{corollary}}
\allowdisplaybreaks[1]

\begin{document}

\begin{abstract}
We give a new short proof that the wheeled operad of unimodular Lie algebras is Koszul and use this to explicitly construct its minimal resolution. A representation of this resolution in a finite dimensional vector space $V$ we call a unimodular $L_\infty$-algebra. Such a structure corresponds to a homological vector field on $V$ together with an invariant measure. We present explicit formulae for homotopy transferred structures, define the deformation complex and give a cohomological obstruction to the extension of an arbitrary structure of finite dimensional $L_\infty$-algebra to a structure of unimodular $L_\infty$-algebra.
\end{abstract}

\maketitle

\section{Introduction}

The theory of operads and props is an effective tool for constructing homotopy theories of many algebraic structures. Given a particular structure (e.g.~associative algebras, Lie algebras, bialgebras), the first step in the construction would be to find an operad or prop $\cP$ whose representations are in one-to-one correspondence with all possible structures of this type. The second step would then be to construct a minimal resolution $\cP_\infty$ of $\cP$. This is often difficult to carry through unless $\cP$ is Koszul, in which case $\cP_\infty$ is the cobar construction on the Koszul dual of $\cP$. The philosophy is that a generic representation of $\cP_\infty$ is the correct up-to-homotopy version of the structure one began with. This program has successfully been applied to many types of algebraic structures, for example, there is an operad $\Lie$ whose representations are precisely Lie algebras. Moreover, $\Lie$ is Koszul so it is possible to compute its minimal resolution $\Lie_\infty$ explicitly. When considering representations of $\Lie_\infty$ one gets $L_\infty$-algebras, a notion which by now has a well-established place in mathematics and mathematical physics. The aim of this paper is to construct the homotopy theory of \emph{unimodular} Lie algebras according to the program just described. The unimodularity condition implies that this can not be done within the category of ordinary operads or props. We show, however, that in the category of wheeled operads, introduced in \cite{Merkulov2007}, unimodularity can be appropriatly handled.

The geometric interpretation of an $L_\infty$-algebra structure on a vector space $V$ is as a homological vector field $Q$ in a formal neighbourhood of $0$ in $V$ which vanishes at $0$. We show that, with this perspective, a structure of unimodular $L_\infty$-algebra is equivalent to a pair $(Q,\rho_f)$ where $\rho_f$ is a $Q$-invariant measure. Such pairs have recently been obtained by P.~Mnëv in \cite{Mnev2006} via Batalin-Vilkovisky quantization of a certain extended BF-theory. Our results imply that the Feynman integrals in \cite{Mnev2006} assemble precisely into a representation of a \emph{minimal} resolution of the wheeled operad of unimodular Lie algebras, and hence indeed have homotopy theoretical meaning.

A few words on notation and conventions. Our ground field is denoted $\k$ and is of characteristic zero. The set $\{1,\ldots,n\}$ is denoted $[n]$ and $\Sigma_n$ is the group of bijections of $[n]$. The trivial representation of $\Sigma_n$ is denoted $\one_n$, while the sign representation is denoted $\sgn_n$. For a permutation $\sigma$, $\sgn\sigma$, and also $(-1)^\sigma$, denote the sign of the permutation. For $[p_1+\cdots+p_n]=I_1\sqcup\cdots\sqcup I_n$, with $|I_j|=\{i_{j,1}<\cdots<i_{j,p_j}\}$,
\[
\sigma_{I_1,\ldots,I_n}=
\begin{pmatrix}
1 & \ldots & p_1 & p_1+1 & \ldots & p_1+p_2 & \ldots & p_1+\cdots+p_n \\
i_{1,1} & \ldots & i_{1,p_1} & i_{2,1} & \ldots & i_{2,p_2} & \ldots & i_{n,p_n}
\end{pmatrix}.
\]
The abbreviation dg stands for differential graded and all differentials are of degree $+1$. The term `corolla' refers to a graph with a single vertex.

The paper is organized as follows. In Section $2$ we recall the theory of wheeled properads and operads, the bar and cobar constructions and Koszulness. In Section $3$ we introduce the wheeled operad of unimodular Lie algebras $\ULie$ and show that it is Koszul. This gives us explicitly its minimal resolution $\ULie_\infty$ in terms of generators and differential. In Section $4$ we give explicit formulae for the transfered structure of unimodular $L_\infty$-algebra. In Section $5$ we interpret a structure of unimodular $L_\infty$-algebra as a homological vector field together with an invariant measure. In Section $6$ we define the deformation complex of a unimodular $L_\infty$-algebra $\fg$. It is a dg Lie algebra and we show that its Maurer-Cartan elements correspond deformations of $\fg$. In Section $7$ we define the characteristic class of an $L_\infty$-algebra and show that we can extend the structure to a unimodular $L_\infty$-algebra if and only if the characteristic class vanishes. Finally, for the sake of completeness, we address the topic of characteristic classes of unimodular associative algebras in an Appendix.
\section{Recollection of the theory of wheeled operads}

\subsection{Decorated graphs and wheeled operads}\hspace{1pt}

Wheeled operads are particular cases of wheeled properads. We find it easier to define this more general notion first and then focus on the special case at hand. For reference on properads and their Koszulness, see \cite{Vallette2007}. For their wheeled extensions, see \cite{Markl2007}.

We denote with $\fG_c$ the class of directed, connected, labeled graphs and by $\fG_c(m,n)$ the subclass of graphs with $m$ outgoing, and $n$ ingoing, legs. By `labeled' we mean that $G\in\fG_c(m,n)$ is equipped with bijections $out(G)\ra[m]$ and $in(G)\ra[n]$, where $out(G)$ and $in(G)$ denote the set of outgoing respectively ingoing legs. We draw graphs directed from top to bottom.

\bex
\[
\begin{xy}
(0,0);(-4,4)**@{-},
(0,0);(4,4)**@{-},
(0,0);(0,-4)**@{-},
(-4,6)*{\scriptstyle{1}},
(4,6)*{\scriptstyle{2}},
(0,-6)*{\scriptstyle{1}}
\end{xy}
\in\fG_c(1,2)\text{ and }
\begin{xy}
(-6,-6);(2,2)**@{-},
(2,-6);(-6,2)**@{-},
(2,2);(6,-2)**@{-},
(2,2);(2,6)**@{-},
(2,6);(2,-6)**\crv{(2,10)&(10,0)&(6,-10)},
(-6,4)*{\scriptstyle{1}},
(6,-4)*{\scriptstyle{1}},
(-6,-8)*{\scriptstyle{2}}
\end{xy}
\in\fG_c(2,1).
\]
\eex

For a graph $G\in\fG_c$, we denote its set of vertices with $\mathbf{v}(G)$. Given a vertex $v$, its set of outgoing legs and edges is denoted $out(v)$ and its set of ingoing legs and edges is denoted $in(v)$. The biarity of a vertex $v$ is $(|out(v)|,|in(v)|)$. We have $\mathbf{v}(G)=\mathbf{v}_c(G)\sqcup\mathbf{v}_{nc}(G)$, where $\mathbf{v}_c(G)$ is the set of cyclic vertices, i.e.~those which lie on a closed directed path, and $\mathbf{v}_{nc}(G)$ is the set of non-cyclic vertices. The edges of $G$ is denoted $\mathbf{e}(G)$ and, as for vertices, $\mathbf{e}(G)=\mathbf{e}_c(G)\sqcup\mathbf{e}_{nc}(G)$.

Recall that a dg $\Sigma$-bimodule $\cE$ is a collection $\{\cE(m,n)\}$, for $m,n\geq 0$, of dg vector spaces with compatible left $\Sigma_m$- and right $\Sigma_n$-actions.

Given a $k$-vertex graph $G$ and a $\Sigma$-bimodule $\cE$, let $v\in\mathbf{v}(G)$ with $|out(v)|=p$ and $|in(v)|=q$. The vector space $\k[\Bij([p],out(v))]$ spanned by all bijections $[p]\ra out(v)$ is naturally a right $\Sigma_p$-module. Like-wise $\k[\Bij(in(v),[q])]$ is naturally a left $\Sigma_q$-module. We define the space of \emph{decorations of the vertex $v$ by $\cE$} as
$$
\cE(out(v),in(v))=\k[\Bij([p],out(v))]\otimes_{\Sigma_p}\cE(p,q)\otimes_{\Sigma_q}\k[\Bij(in(v),[q])].
$$
The space of \emph{decorations of the graph $G$ by $\cE$} is defined as
$$
\decor{G}{\cE}=\left(\bigotimes_{v\in\mathbf{v}(G)}\cE(out(v),in(v))\right)_{\Aut G},
$$
where
$$
\bigotimes_{v\in\mathbf{v}(G)}\cE(out(v),in(v))=\left(\bigoplus_{\gamma\in\Bij([k],\mathbf{v}(G))}\bigotimes_{i=1}^k\cE(out(\gamma(i)),in(\gamma(i)))\right)_{\Sigma_k}
$$
is the unordered tensor product and $\Aut G$ is the group of automorphisms of $G$ which fix legs. For $G\in\fG_c(m,n)$, relabeling gives $\decor{G}{\cE}$ a left $\Sigma_m$- and a right $\Sigma_n$-action. An element of $\decor{G}{\cE}$, for a $k$-vertex graph $G$ is essentially a pair $(G,[e_1\otimes\cdots\otimes e_k])$ of the graph and an equivalence class of a tensor product of elements of $\cE$.

For any connected subgraph $H\subset G$, we denote with $G/H$ the graph in which all vertices and edges of $H$ has been contracted into a single vertex. Given a $\Sigma$-bimodule $\cE$, let $\{\mu_G:\decor{G}{\cE}\ra\cE\}_{G\in\fG_c}$ be a collection of equivatiant linear maps. For a connected $H$ there is an induced map $\mu^G_H:\decor{G}{\cE}\ra\decor{G/H}{\cE}$, which equals $\mu_H$ on the vertices of $H$ and the identity elsewhere.

\bde
A \emph{wheeled properad} is a $\Sigma$-bimodule $\cP$ together with a collection of equivariant linear maps $\{\mu_G:\decor{G}{\cP}\ra\cP\}_{G\in\fG_c}$ such that, for any connected subgraph $H\subset G$, $\mu_G=\mu_{G/H}\circ\mu^G_H$. A \emph{wheeled operad} is a wheeled properad $\cP$ such that if $m\geq 2$, then $\cP(m,n)=0$.
\ede

From the equations satisfied by the collection $\{\mu_G\}$ it follows that it is determined by the subcollection of those maps $\mu_G$ for which $G$ is either a two vertex graph without directed cycles or a corolla with a loop.

Given a wheeled operad $\cP$ its \emph{operadic part}, $\cP_o$, is defined by
$$
\cP_o(n)=\cP(1,n).
$$
It is clear that the $\Sigma$-module $\{\cP_o(n)\}$ is an ordinary operad. The \emph{wheeled part} of $\cP$, $\cP_w$, is defined by
$$
\cP_w(n)=\cP(0,n)
$$
and is, obviously, a right $\cP_o$-module.

The inclusion of genus zero graphs into $\fG_c$ induces a forgetful functor from wheeled properads to ordinary properads. This functor has an adjoint $\cP\mapsto\cP^\cal$, called \emph{wheelification}. The wheelification of an ordinary operad is a wheeled operad.

Dually, consider a dg $\Sigma$-bimodule $\cE$ and a collection of equivariant linear maps $\{\Delta_G:\cE\ra\decor{G}{\cE}\}_{G\in\fG_c}$. For any connected subgraph $H\subset G$ there is an induced map $\Delta^G_H:\decor{G/H}{\cE}\ra\decor{G}{\cE}$, which equals $\Delta_H$ on the vertex of $G/H$ corresponding to $H$ and the identity elsewhere.

\bde
A \emph{wheeled coproperad} is a $\Sigma$-bimodule $\cC$ together with a collection of linear maps $\{\Delta_G:\cC\ra\decor{G}{\cC}\}_{G\in\fG_c}$ such that, for any connected subgraph $H\subset G$, $\Delta_G=\Delta^G_H\circ\Delta_{G/H}$. A \emph{wheeled cooperad} is a wheeled coproperad $\cC$ such that if $m\geq 2$, then $\cP(m,n)=0$.
\ede

Again, $\{\Delta_G\}$ is determined by the subcollections of maps $\Delta_G$ for which $G$ is either a two vertex graph without directed cycles or a corolla with a loop.

\bex
Given a dg $\Sigma$-bimodule $\cE$, the \emph{free wheeled properad} on $\cE$, $\cF^\cal(\cE)$, has as underlying dg $\Sigma$-bimodule
$$
\cF^\cal(m,n)=\bigoplus_{G\in\fG_c(m,n)}\decor{G}{\cE},
$$
with differential $\partial_\cE$ induced by that of $\cE$. An element of $\decor{G}{\cF^\cal(\cE)}$ is a graph with vertices decorated by graphs decorated by $\cE$. The sturcture of wheeled properad on $\cF^\cal(\cE)$ forgets this double decoration. Note that $\fG_c$ contains the two exceptional graphs $\downarrow$ and $\cal$ without vertices.

The \emph{cofree wheeled coproperad} on $\cE$, $\cF^\cal_c(\cE)$, has the same underlying dg $\Sigma$-bimodule as $\cF^\cal(\cE)$. The morphism $\Delta_G$ on a corolla decorated by some graph $G'$ is the sum of decorations of $G$ by connected subgraphs of $G'$ such that when forgetting the double decoration one gets $G'$.
\eex

Note that if $m\geq 2$ implies $\cE(m,n)=0$, then $\cF^\cal(\cE)$ is a wheeled operad and $\cF^\cal_c(\cE)$ is a wheeled cooperad.

\bex
Given a finite dimensional dg vector space $(V,d)$, let $\oEnd^\cal_V$ be defined by
$$
\oEnd^\cal_V(m,n)=
\begin{cases}
\Hom(V^{\otimes n},\k) & m=0, \\
\Hom(V^{\otimes n},V) & m=1, \\
0 & \text{otherwise}.
\end{cases}
$$
It is a dg $\Sigma$-bimodule with the usual grading, differential and $\Sigma$-action. For any $f:V^{\otimes n+1}\ra V$ and $1\leq i\leq n+1$ we may fix elements $x_1,\ldots,x_{i-1},x_{i+1},\ldots,x_{n+1}\in V$ and get
$$
f(x_1,\ldots,x_{i-1},-,x_{i+1},\ldots,x_{n+1}):V\ra V.
$$
Taking the trace of this endomorphism gives us a map $V^{\otimes n}\ra\k$. This, together with usual composition of homomorphisms, gives $\oEnd^\cal_V$ the structure of dg wheeled operad. It is called the \emph{wheeled endomorphism operad} of $V$.
\eex

\bde
A \emph{representation} of a dg wheeled operad $\cP$ in a finite dimensional dg vector space $V$ is a morphism of dg wheeled operads $\cP\ra\oEnd^\cal_V$. We say that a representation gives $V$ the structure of a \emph{$\cP$-algebra}.
\ede

For a graded vector space $V=\oplus_i V^i$, $V[j]=\oplus_i V[j]^i$ denotes the \emph{shifted} vector space with $V[j]^i=V^{j+i}$. On the category of wheeled operads there is a \emph{parity change functor}, $\cP\mapsto\Pi\cP$, such that representations of $\cP$ in a vector space $V$ are in one-to-one correspondence with representations of $\Pi\cP$ in $V[1]$. Explicitly,
$$
\Pi\cP(m,n)=\sgn_m\otimes\cP(m,n)\otimes\sgn_n[m-n],
$$
so canonically $\Pi\oEnd_V^\cal=\oEnd_{V[1]}^\cal$ by sending $f:V^{\otimes n}\ra V^{\otimes m}$ to
$$
(\underbrace{s^{-1}\otimes\cdots\otimes s^{-1}}_m)f(\underbrace{s\otimes\cdots\otimes s}_n):V[1]^{\otimes n}\ra V[1]^{\otimes m}.
$$
Here $s$ denotes the degree $1$ map $V[1]\ra V$ identifying $V[1]^i$ with $V^{i+1}$.

\subsection{Bar and cobar construction for wheeled operads}\hspace{1pt}

Consider the wheeled operad
$$
I^\circlearrowleft(m,n)=
\begin{cases}
\k & \text{if }(m,n)=(0,0)\text{ or }(m,n)=(1,1), \\
0 & \text{otherwise}.
\end{cases}
$$

\bde
An \emph{augmentation} of a wheeled operad $\cP$ is a morphism $\cP\ra I^\circlearrowleft$. An \emph{augmented wheeled operad} is a wheeled operad together with an augmentation. Given an augmented wheeled operad $\cP$, the kernel of the augmentation is denoted $\bar{\cP}$ and is called the \emph{augmentation ideal}.
\ede

\bde
The \emph{wheeled suspension} of a dg $\Sigma$-bimodule $\cE$, denoted $\w\cE$, is defined as
$$
\w\cE(m,n)=\cE(m,n)[2m-n]\otimes\sgn_n
$$
and the \emph{wheeled desuspension}, $\w^{-1}\cE$, as
$$
\w^{-1}\cE(m,n)=\cE(m,n)[n-2m]\otimes\sgn_n.
$$
\ede

Given a dg wheeled properad $\cP$, the collection of maps $\mu_G$ for which $G$ is either a two vertex graph without directed cycles or a corolla with a loop defines a morphism $\mu:\cF^\cal(\cP)\ra\cP$. This $\mu$ induces a degree one morphism $\cF^\cal_c(\mathbf{w}\bar{\cP})\ra\mathbf{w}\bar{\cP}$, which in turn determines a degree one coderivation $\partial_\mu$ of $\cF^\cal_c(\mathbf{w}\bar{\cP})$. One can check that the equations satisfied by the $\mu_G$ are equivalent with $\partial_\mu$ being a codifferential.

\bde
The \emph{wheeled bar construction} on a dg wheeled properad $\cP$ is defined as
$$
(B^\circlearrowleft(\cP),\partial_B)=(\cF_c^\circlearrowleft(\w\bar{\cP}),\partial_\cP+\partial_\mu).
$$
\ede

Dually, given a dg wheeled coproperad $\cC$, the collection $\Delta_G$ for which $G$ is either a two vertex graph without directed cycles or a corolla with a loop defines a morphism $\Delta:\cC\ra\cF^\cal_c(\cC)$. This induces a degree one morphism $\mathbf{w}^{-1}\cC\ra\cF^\cal(\mathbf{w}^{-1}\cC)$, which determines a degree one derivation $\partial_\Delta$ of $\cF^\cal(\mathbf{w}^{-1}\cC)$. One checks that the equations satisfied by the $\Delta_G$ are equivalent with $\partial_\Delta$ being a differential.

\bde
The \emph{wheeled cobar construction} on a dg wheeled coproperad $\cC$ is defined as
$$
(\Omega^\circlearrowleft(\cC),\partial_\Omega)=(\cF^\circlearrowleft(\w^{-1}\bar{\cC}),\partial_\cC+\partial_\Delta).
$$
\ede

\subsection{Quadratic wheeled operads and wheeled Koszulness}\hspace{1pt}

\bde
A wheeled operad $\cP$ is \emph{quadratic} if it has a presentation
$$
\cP=\cF^\circlearrowleft(\cE)/(R),
$$
where $\cE$ is non-zero only for $(m,n)=(1,2)$ and $R$ is a subspace of $\cF^\circlearrowleft(\cE)(1,3)\oplus\cF^\circlearrowleft(\cE)(0,1)$.
\ede

Note that any quadratic wheeled operad is naturally augmented by projection onto the part spanned by the exceptional graphs.

\bex
Let $\cE(1,2)=\one_2$ and denote its generator by
$$
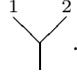
\xygraph{
!{<0pt,0pt>;<10pt,0pt>:<0pt,10pt>::},
!{(0,0)}-!{(0,-1)},
!{(0,0)}-!{(-1,1)},
!{(0,0)}-!{(1,1)},
!{(-1,1.4)}*{\scriptstyle{1}},
!{(1,1.4)}*{\scriptstyle{2}}
}.
$$
Then the wheelification of the operad $\Com$ of commutative associative algebras has the presentation
$$
\Com^\cal=
\frac{
\cF^\circlearrowleft\left(
\xygraph{
!{<0pt,0pt>;<10pt,0pt>:<0pt,10pt>::},
!{(0,0)}-!{(0,-1)},
!{(0,0)}-!{(-1,1)},
!{(0,0)}-!{(1,1)},
!{(-1,1.4)}*{\scriptstyle{1}},
!{(1,1.4)}*{\scriptstyle{2}}
}
\right)
}{
\left(
\xygraph{
!{<0pt,0pt>;<10pt,0pt>:<0pt,10pt>::},
!{(-0.5,-0.5)}-!{(-0.5,-1.5)},
!{(-0.5,-0.5)}-!{(-1.5,0.5)},
!{(-0.5,-0.5)}-!{(1.5,1.5)},
!{(0.5,0.5)}-!{(-0.5,1.5)},
!{(-1.5,0.9)}*{\scriptstyle{3}},
!{(-0.5,1.9)}*{\scriptstyle{1}},
!{(1.5,1.9)}*{\scriptstyle{2}}
}-
\xygraph{
!{<0pt,0pt>;<10pt,0pt>:<0pt,10pt>::},
!{(-0.5,-0.5)}-!{(-0.5,-1.5)},
!{(-0.5,-0.5)}-!{(-1.5,0.5)},
!{(-0.5,-0.5)}-!{(1.5,1.5)},
!{(0.5,0.5)}-!{(-0.5,1.5)},
!{(-1.5,0.9)}*{\scriptstyle{1}},
!{(-0.5,1.9)}*{\scriptstyle{2}},
!{(1.5,1.9)}*{\scriptstyle{3}}
},
\xygraph{
!{<0pt,0pt>;<10pt,0pt>:<0pt,10pt>::},
!{(-0.5,-0.5)}-!{(-0.5,-1.5)},
!{(-0.5,-0.5)}-!{(-1.5,0.5)},
!{(-0.5,-0.5)}-!{(1.5,1.5)},
!{(0.5,0.5)}-!{(-0.5,1.5)},
!{(-1.5,0.9)}*{\scriptstyle{2}},
!{(-0.5,1.9)}*{\scriptstyle{1}},
!{(1.5,1.9)}*{\scriptstyle{3}}
}-
\xygraph{
!{<0pt,0pt>;<10pt,0pt>:<0pt,10pt>::},
!{(-0.5,-0.5)}-!{(-0.5,-1.5)},
!{(-0.5,-0.5)}-!{(-1.5,0.5)},
!{(-0.5,-0.5)}-!{(1.5,1.5)},
!{(0.5,0.5)}-!{(-0.5,1.5)},
!{(-1.5,0.9)}*{\scriptstyle{1}},
!{(-0.5,1.9)}*{\scriptstyle{3}},
!{(1.5,1.9)}*{\scriptstyle{2}}
}
\right)
}.
$$
\eex

\bde
The \emph{\v Cech dual} of a $\Sigma$-bimodule $\cE$, denoted $\check{\cE}$, is defined by
$$
\check{\cE}(m,n)=\sgn_m\otimes \cE(m,n)^*\otimes\sgn_n.
$$
\ede

\bde
The \emph{wheeled Kozul dual operad} of the quadratic wheeled operad $\cP=\cF^\circlearrowleft(\cE)/(R)$, is defined as
$$
\cP^!=\cF^\circlearrowleft(\check{\cE})/(R^\perp).
$$
Here $R^\perp$ is the kernel of the composition
$$
\cF^\circlearrowleft(\check{\cE})(1,3)\oplus\cF^\circlearrowleft(\check{\cE})(0,1)\isoto(\cF^\circlearrowleft(\cE)(1,3)\oplus\cF^\circlearrowleft(\cE)(0,1))^*\ra R^*.
$$
\ede

The operadic part of the wheeled Koszul dual operad of a wheeled operad $\cP$ is the ordinary Koszul dual operad of the operadic part of $\cP$, i.e.
$$
(\cP^!)_o=(\cP_o)^!.
$$
The wheeled part equals the quotient of the wheeled part of the wheeled completion of $\cP^!_o$ by the ideal generated by $R^\perp_w$,
$$
(\cP^!)_w=((\cP^!_o)^\circlearrowleft)_w/(R^\perp_w).
$$

\bde
The \emph{wheeled Koszul dual cooperad} of the quadratic wheeled operad $\cP$ is defined as
$$
\cP^\text{!`}=H^0(B^\circlearrowleft(\cP),\partial_B).
$$
\ede
One can show that $\cP^\text{!`}$ is a wheeled subcooperad of $B^\cal(\cP)$ and that $\cP^!=(\cP^\text{!`})^*$.

\bde
A quadratic wheeled operad $\cP$ is \emph{wheeled Koszul} if the inclusion
$$
(\cP^\text{!`},0)\ra(B^\circlearrowleft(\cP),\partial_B)
$$
is a quasi-isomorphism, i.e.~induces an isomorphism on cohomology.
\ede

The wheeled operad $\cP$ is wheeled Koszul if and only if $\cP^!$ is wheeled Koszul. Moreover, $\cP$ is wheeled Koszul if and only if the natural projection
$$
(\Omega^\circlearrowleft(\cP^\text{!`}),\partial_\Omega)\ra(\cP,0)
$$
is a quasi-isomorphism. Hence, if $\cP$ is Koszul, then $(\Omega^\circlearrowleft(\cP^\text{!`}),\partial_\Omega)$ is a free minimial resolution of $\cP$. In this case we write $\cP_\infty=(\Omega^\circlearrowleft(\cP^\text{!`}),\partial_\Omega)$.

\bde
An ordinary Koszul operad $\cP$ is \emph{stably Koszul} if the wheelification of the quasi-isomorphism $(\Omega(\cP^\text{!`}),\partial_\Omega)\ra(\cP,0)$ is still a quasi-isomorphism.
\ede

In \cite{Merkulov2007} it is shown that $\Lie$ is stably Koszul. We shall use this fact in the proof of wheeled Koszulness for the wheeled operad of unimodular Lie algebras, $\ULie$.

\section{Unimodular Lie algebras}

Any Lie algebra $(\fg,[-,-])$ is endowed with a canonical representation in itself, the \emph{adjoint representation} $\ad:\fg\ra\End(\fg)$, defined by $\ad(x)(y)=[x,y]$. A finite dimensional Lie algebra $\fg$ is \emph{unimodular} if $\tr(\ad x)=0$ for all $x\in\fg$. Important examples of unimodular Lie algebras include semi-simple Lie algebras. If we choose a basis $\{e_\alpha\}$ for $\fg$, the Lie algebra structure is determined by the structure constants $L^\gamma_{\alpha\beta}$ defined by
$$
[e_\alpha,e_\beta]=L^\gamma_{\alpha\beta}e_\gamma.
$$
Unimodularity is then expressed as $L^\beta_{\alpha\beta}=0$ for all $\alpha$.

Unimodular Lie algebras are representations of a quadratic wheeled operad, $\ULie=\cF^\circlearrowleft(\cE)/(R)$, where $\cE(1,2)=\sgn_2$ and $R$ is spanned by the Jacobi relation and the trace condition. We denote the generator of $\cE(1,2)$ by
$$
\xygraph{
!{<0pt,0pt>;<10pt,0pt>:<0pt,10pt>::},
!{(0,0)}-!{(0,-1)},
!{(0,0)}-!{(-1,1)},
!{(0,0)}-!{(1,1)},
!{(-1,1.4)}*{\scriptstyle{1}},
!{(1,1.4)}*{\scriptstyle{2}}
}
$$
so that
$$
\ULie=
\frac{
\cF^\circlearrowleft\left(
\xygraph{
!{<0pt,0pt>;<10pt,0pt>:<0pt,10pt>::},
!{(0,0)}-!{(0,-1)},
!{(0,0)}-!{(-1,1)},
!{(0,0)}-!{(1,1)},
!{(-1,1.4)}*{\scriptstyle{1}},
!{(1,1.4)}*{\scriptstyle{2}}
}
\right)
}{
\left(
\xygraph{
!{<0pt,0pt>;<10pt,0pt>:<0pt,10pt>::},
!{(-0.5,-0.5)}-!{(-0.5,-1.5)},
!{(-0.5,-0.5)}-!{(-1.5,0.5)},
!{(-0.5,-0.5)}-!{(1.5,1.5)},
!{(0.5,0.5)}-!{(-0.5,1.5)},
!{(-1.5,0.9)}*{\scriptstyle{1}},
!{(-0.5,1.9)}*{\scriptstyle{2}},
!{(1.5,1.9)}*{\scriptstyle{3}}
}-
\xygraph{
!{<0pt,0pt>;<10pt,0pt>:<0pt,10pt>::},
!{(-0.5,-0.5)}-!{(-0.5,-1.5)},
!{(-0.5,-0.5)}-!{(-1.5,0.5)},
!{(-0.5,-0.5)}-!{(1.5,1.5)},
!{(0.5,0.5)}-!{(-0.5,1.5)},
!{(-1.5,0.9)}*{\scriptstyle{2}},
!{(-0.5,1.9)}*{\scriptstyle{1}},
!{(1.5,1.9)}*{\scriptstyle{3}}
}+
\xygraph{
!{<0pt,0pt>;<10pt,0pt>:<0pt,10pt>::},
!{(-0.5,-0.5)}-!{(-0.5,-1.5)},
!{(-0.5,-0.5)}-!{(-1.5,0.5)},
!{(-0.5,-0.5)}-!{(1.5,1.5)},
!{(0.5,0.5)}-!{(-0.5,1.5)},
!{(-1.5,0.9)}*{\scriptstyle{3}},
!{(-0.5,1.9)}*{\scriptstyle{1}},
!{(1.5,1.9)}*{\scriptstyle{2}}
},
\xygraph{
!{<0pt,0pt>;<10pt,0pt>:<0pt,10pt>::},
!{(0,0)}-!{(0,-1)},
!{(0,0)}-!{(-1,1)},
!{(0,0)}-!{(1,1)},
!{(0,-1)}-@(d,ur)!{(1,1)}
}
\right)
}.
$$

We will show that $\ULie$ is wheeled Koszul by determining $(\Omega^\circlearrowleft(\ULie^\text{!`}),\partial)$ and showing that its cohomology equals $\ULie$. Since $\ULie^\text{!`}=(\ULie^!)^*$, our first step will be to determine $\ULie^!$.

\ble
The wheeled Koszul dual operad of $\ULie$ is $\Com^\circlearrowleft$.
\ele
\bpro\label{prop:threeone}
We have $\ULie_o=\Lie$ so that $(\ULie^!)_o=\Lie^!=\Com$. Next, since $R_w$ spans all of $\cF^\circlearrowleft(\cE)(0,1)$, $R_w^\perp=0$ and $(\ULie^!)_w=\Com^\circlearrowleft_w$. Hence, $\ULie^!=\Com^\circlearrowleft$.
\epro

Since $\ULie_o=\Lie$ and $\Lie$ is Koszul, we have that $\Omega^\cal(\ULie^\text{!`})_o=\Lie_\infty$. Hence, only the wheeled part of $\Omega^\cal(\ULie^\text{!`})$ will concern us. Since $\Com^\cal_w(n)=\one_n$ for $n\geq 1$, $\Omega^\circlearrowleft(\ULie^\text{!`})_w$ has a single skew-symmetric generator in each arity $n\geq 1$. We denote it
$$
\xygraph{
!{<0pt,0pt>;<10pt,0pt>:<0pt,10pt>::},
!{(0,0)}-!{(-1,1)},
!{(0,0)}-!{(-0.5,1)},
!{(0,0)}-!{(0,1)},
!{(0,0)}-!{(0.5,1)},
!{(0,0)}-!{(1,1)},
!{(-1,1.4)}*{\scriptstyle{1}},
!{(1,1.4)}*{\scriptstyle{n}},
!{(0,0)}*{-}
}.
$$
Considering the cooperad structure in $(\Com^\circlearrowleft)^*$, the differential $\partial$ applied to this generator has terms of the form
$$
\xygraph{
!{<0pt,0pt>;<10pt,0pt>:<0pt,10pt>::},
!{(-0.5,-0.5)}-!{(-1.5,0.5)},
!{(-0.5,-0.5)}-!{(-1,0.5)},
!{(-0.5,-0.5)}-!{(-0.5,0.5)},
!{(-0.5,-0.5)}-!{(0,0.5)},
!{(-0.5,-0.5)}-!{(1.5,1.5)},
!{(0.5,0.5)}-!{(-0.5,1.5)},
!{(0.5,0.5)}-!{(0,1.5)},
!{(0.5,0.5)}-!{(0.5,1.5)},
!{(0.5,0.5)}-!{(1,1.5)},
!{(-1,0.9)}*{\scriptstyle{I}},
!{(0.5,1.9)}*{\scriptstyle{J}},
!{(-0.5,-0.5)}*{-}
}
\text{ where }I\sqcup J=[n]\text{ and }0\leq |I|\leq n-2,\text{ and }
\xygraph{
!{<0pt,0pt>;<10pt,0pt>:<0pt,10pt>::},
!{(0,0)}-!{(0,-1)},
!{(0,0)}-!{(-1,1)},
!{(0,0)}-!{(-0.5,1)},
!{(0,0)}-!{(0,1)},
!{(0,0)}-!{(0.5,1)},
!{(0,0)}-!{(1,1)},
!{(0,-1)}-@(d,ur)!{(1,1)},
!{(-1,1.4)}*{\scriptstyle{1}},
!{(0.5,1.4)}*{\scriptstyle{n}}
}.
$$
Hence, to determine $(\Omega^\circlearrowleft(\ULie^\text{!`}),\partial)$, we need to find signs in front of these terms such that $\partial^2=0$.

\bpr
$(\Omega^\circlearrowleft(\ULie^\text{!`}),\partial)$ is free on generators
$$
\left\{
\xygraph{
!{<0pt,0pt>;<10pt,0pt>:<0pt,10pt>::},
!{(0,0)}-!{(0,-1)},
!{(0,0)}-!{(-1,1)},
!{(0,0)}-!{(-0.5,1)},
!{(0,0)}-!{(0,1)},
!{(0,0)}-!{(0.5,1)},
!{(0,0)}-!{(1,1)},
!{(-1,1.4)}*{\scriptstyle{1}},
!{(1,1.4)}*{\scriptstyle{n}}
}
\right\}_{n=2}^\infty
\text{ and }
\left\{
\xygraph{
!{<0pt,0pt>;<10pt,0pt>:<0pt,10pt>::},
!{(0,0)}-!{(-1,1)},
!{(0,0)}-!{(-0.5,1)},
!{(0,0)}-!{(0,1)},
!{(0,0)}-!{(0.5,1)},
!{(0,0)}-!{(1,1)},
!{(-1,1.4)}*{\scriptstyle{1}},
!{(1,1.4)}*{\scriptstyle{n}},
!{(0,0)}*{-}
}
\right\}_{n=1}^\infty,
$$
skew-symmetric and of degrees $2-n$ and $-n$ respectively. The differential is defined by
\begin{align*}
\partial
\xygraph{
!{<0pt,0pt>;<10pt,0pt>:<0pt,10pt>::},
!{(0,0)}-!{(0,-1)},
!{(0,0)}-!{(-1,1)},
!{(0,0)}-!{(-0.5,1)},
!{(0,0)}-!{(0,1)},
!{(0,0)}-!{(0.5,1)},
!{(0,0)}-!{(1,1)},
!{(-1,1.4)}*{\scriptstyle{1}},
!{(1,1.4)}*{\scriptstyle{n}}
}=&
\sum_{\substack{I\sqcup J=[n] \\ 1\leq |I|\leq n-2}}
(-1)^{\sigma_{I,J}+|I||J|}
\xygraph{
!{<0pt,0pt>;<10pt,0pt>:<0pt,10pt>::},
!{(-0.5,-0.5)}-!{(-0.5,-1.5)},
!{(-0.5,-0.5)}-!{(-1.5,0.5)},
!{(-0.5,-0.5)}-!{(-1,0.5)},
!{(-0.5,-0.5)}-!{(-0.5,0.5)},
!{(-0.5,-0.5)}-!{(0,0.5)},
!{(-0.5,-0.5)}-!{(1.5,1.5)},
!{(0.5,0.5)}-!{(-0.5,1.5)},
!{(0.5,0.5)}-!{(0,1.5)},
!{(0.5,0.5)}-!{(0.5,1.5)},
!{(0.5,0.5)}-!{(1,1.5)},
!{(-1,0.9)}*{\scriptstyle{I}},
!{(0.5,1.9)}*{\scriptstyle{J}}
}, \\
\partial
\xygraph{
!{<0pt,0pt>;<10pt,0pt>:<0pt,10pt>::},
!{(0,0)}-!{(-1,1)},
!{(0,0)}-!{(-0.5,1)},
!{(0,0)}-!{(0,1)},
!{(0,0)}-!{(0.5,1)},
!{(0,0)}-!{(1,1)},
!{(-1,1.4)}*{\scriptstyle{1}},
!{(1,1.4)}*{\scriptstyle{n}},
!{(0,0)}*{-}
}=&
\sum_{\substack{I\sqcup J=[n] \\ 0\leq |I|\leq n-2}}
(-1)^{\sigma_{I,J}+|I||J|}
\xygraph{
!{<0pt,0pt>;<10pt,0pt>:<0pt,10pt>::},
!{(-0.5,-0.5)}-!{(-1.5,0.5)},
!{(-0.5,-0.5)}-!{(-1,0.5)},
!{(-0.5,-0.5)}-!{(-0.5,0.5)},
!{(-0.5,-0.5)}-!{(0,0.5)},
!{(-0.5,-0.5)}-!{(1.5,1.5)},
!{(0.5,0.5)}-!{(-0.5,1.5)},
!{(0.5,0.5)}-!{(0,1.5)},
!{(0.5,0.5)}-!{(0.5,1.5)},
!{(0.5,0.5)}-!{(1,1.5)},
!{(-1,0.9)}*{\scriptstyle{I}},
!{(0.5,1.9)}*{\scriptstyle{J}},
!{(-0.5,-0.5)}*{-}
}+
\xygraph{
!{<0pt,0pt>;<10pt,0pt>:<0pt,10pt>::},
!{(0,0)}-!{(0,-1)},
!{(0,0)}-!{(-1,1)},
!{(0,0)}-!{(-0.5,1)},
!{(0,0)}-!{(0,1)},
!{(0,0)}-!{(0.5,1)},
!{(0,0)}-!{(1,1)},
!{(0,-1)}-@(d,ur)!{(1,1)},
!{(-1,1.4)}*{\scriptstyle{1}},
!{(0.5,1.4)}*{\scriptstyle{n}}
}.
\end{align*}
\epr
\bpro
We must show that $\partial^2=0$. Since $\Omega^\circlearrowleft(\ULie^\text{!`})_o=\Lie_\infty$, only the generators of the second type concern us. Consider first the sum
\begin{align*}
\sum_{I'\sqcup J'=[n+1]}(-1)^{\sigma_{I',J'}+|I'||J'|}
\xygraph{
!{<0pt,0pt>;<10pt,0pt>:<0pt,10pt>::},
!{(-0.5,-0.5)}-!{(-1.5,0.5)},
!{(-0.5,-0.5)}-!{(-1,0.5)},
!{(-0.5,-0.5)}-!{(-0.5,0.5)},
!{(-0.5,-0.5)}-!{(0,0.5)},
!{(-0.5,-0.5)}-!{(1.5,1.5)},
!{(0.5,0.5)}-!{(-0.5,1.5)},
!{(0.5,0.5)}-!{(0,1.5)},
!{(0.5,0.5)}-!{(0.5,1.5)},
!{(0.5,0.5)}-!{(1,1.5)},
!{(-1,0.9)}*{\scriptstyle{I'}},
!{(0.5,1.9)}*{\scriptstyle{J'}},
!{(-0.5,-0.5)}*{-}
}=&\sum_{\substack{I'\sqcup J'=[n+1] \\ n+1\in I'}}(-1)^{\sigma_{I',J'}+|I'||J'|}
\xygraph{
!{<0pt,0pt>;<10pt,0pt>:<0pt,10pt>::},
!{(-0.5,-0.5)}-!{(-1.5,0.5)},
!{(-0.5,-0.5)}-!{(-1,0.5)},
!{(-0.5,-0.5)}-!{(-0.5,0.5)},
!{(-0.5,-0.5)}-!{(0,0.5)},
!{(-0.5,-0.5)}-!{(1.5,1.5)},
!{(0.5,0.5)}-!{(-0.5,1.5)},
!{(0.5,0.5)}-!{(0,1.5)},
!{(0.5,0.5)}-!{(0.5,1.5)},
!{(0.5,0.5)}-!{(1,1.5)},
!{(-1,0.9)}*{\scriptstyle{I'}},
!{(0.5,1.9)}*{\scriptstyle{J'}},
!{(-0.5,-0.5)}*{-}
} \\
&+
\sum_{\substack{I'\sqcup J'=[n+1] \\ n+1\in J'}}(-1)^{\sigma_{I',J'}+|I'||J'|}
\xygraph{
!{<0pt,0pt>;<10pt,0pt>:<0pt,10pt>::},
!{(-0.5,-0.5)}-!{(-1.5,0.5)},
!{(-0.5,-0.5)}-!{(-1,0.5)},
!{(-0.5,-0.5)}-!{(-0.5,0.5)},
!{(-0.5,-0.5)}-!{(0,0.5)},
!{(-0.5,-0.5)}-!{(1.5,1.5)},
!{(0.5,0.5)}-!{(-0.5,1.5)},
!{(0.5,0.5)}-!{(0,1.5)},
!{(0.5,0.5)}-!{(0.5,1.5)},
!{(0.5,0.5)}-!{(1,1.5)},
!{(-1,0.9)}*{\scriptstyle{I'}},
!{(0.5,1.9)}*{\scriptstyle{J'}},
!{(-0.5,-0.5)}*{-}
}.
\end{align*}
If $n+1\in I'$ we let $I=I'\smallsetminus\{n+1\}$, $J=J'$ and note that $I\sqcup J=[n]$ is such that $(-1)^{\sigma_{I',J'}}=(-1)^{\sigma_{I,J}+|J|}$. If instead $n+1\in J'$ we let $I=I'$, $J=J'\smallsetminus\{n+1\}$ and note that $I\sqcup J=[n]$ is such that $(-1)^{\sigma_{I',J'}}=(-1)^{\sigma_{I,J}}$. Hence, the above sum equals
\begin{align*}
\sum_{I\sqcup J=[n]}(-1)^{\sigma_{I,J}+|I||J|+1}
\xygraph{
!{<0pt,0pt>;<10pt,0pt>:<0pt,10pt>::},
!{(-0.5,-0.5)}-!{(-1.5,0.5)},
!{(-0.5,-0.5)}-!{(-1,0.5)},
!{(-0.5,-0.5)}-!{(-0.5,0.5)},
!{(-0.5,-0.5)}-!{(0,0.5)},
!{(-0.5,-0.5)}-!{(0.5,0.5)},
!{(0,0.5)}-!{(-1,1.5)},
!{(0,0.5)}-!{(-0.5,1.5)},
!{(0,0.5)}-!{(0,1.5)},
!{(0,0.5)}-!{(0.5,1.5)},
!{(0,0.5)}-!{(1,1.5)},
!{(-1,0.9)}*{\scriptstyle{I}},
!{(0,1.9)}*{\scriptstyle{J}},
!{(1.5,0.9)}*{\scriptstyle{n+1}},
!{(-0.5,-0.5)}*{-}
}+
(-1)^{\sigma_{I,J}+|I||J|+|I|}
\xygraph{
!{<0pt,0pt>;<10pt,0pt>:<0pt,10pt>::},
!{(-0.5,-0.5)}-!{(-1.5,0.5)},
!{(-0.5,-0.5)}-!{(-1,0.5)},
!{(-0.5,-0.5)}-!{(-0.5,0.5)},
!{(-0.5,-0.5)}-!{(0,0.5)},
!{(-0.5,-0.5)}-!{(1.5,1.5)},
!{(0.5,0.5)}-!{(-0.5,1.5)},
!{(0.5,0.5)}-!{(0,1.5)},
!{(0.5,0.5)}-!{(0.5,1.5)},
!{(0.5,0.5)}-!{(1,1.5)},
!{(-1,0.9)}*{\scriptstyle{I}},
!{(0.5,1.9)}*{\scriptstyle{J}},
!{(2,1.9)}*{\scriptstyle{n+1}},
!{(-0.5,-0.5)}*{-}
}.
\end{align*}
This gives us that
\begin{align*}
\partial
\xygraph{
!{<0pt,0pt>;<10pt,0pt>:<0pt,10pt>::},
!{(-0.5,-0.5)}-!{(-1.5,0.5)},
!{(-0.5,-0.5)}-!{(-1,0.5)},
!{(-0.5,-0.5)}-!{(-0.5,0.5)},
!{(-0.5,-0.5)}-!{(0,0.5)},
!{(-0.5,-0.5)}-!{(1.5,1.5)},
!{(0.5,0.5)}-!{(-0.5,1.5)},
!{(0.5,0.5)}-!{(0,1.5)},
!{(0.5,0.5)}-!{(0.5,1.5)},
!{(0.5,0.5)}-!{(1,1.5)},
!{(-1,0.9)}*{\scriptstyle{I}},
!{(0.5,1.9)}*{\scriptstyle{J}},
!{(-0.5,-0.5)}*{-}
}=&
\sum_{P\sqcup Q=I}(-1)^{\sigma_{P,Q}+|P||Q|+1}
\xygraph{
!{<0pt,0pt>;<10pt,0pt>:<0pt,10pt>::},
!{(-1,-1.5)}-!{(-2,-0.5)},
!{(-1,-1.5)}-!{(-1.5,-0.5)},
!{(-1,-1.5)}-!{(-1,-0.5)},
!{(-1,-1.5)}-!{(0,0.5)},
!{(-1,-1.5)}-!{(2,1.5)},
!{(-0.5,-0.5)}-!{(-1.5,0.5)},
!{(-0.5,-0.5)}-!{(-1,0.5)},
!{(-0.5,-0.5)}-!{(-0.5,0.5)},
!{(-0.5,-0.5)}-!{(0.5,0.5)},
!{(1,0.5)}-!{(0,1.5)},
!{(1,0.5)}-!{(0.5,1.5)},
!{(1,0.5)}-!{(1,1.5)},
!{(1,0.5)}-!{(1.5,1.5)},
!{(-2,0)}*{\scriptstyle{P}},
!{(-0.5,1)}*{\scriptstyle{Q}},
!{(1,2)}*{\scriptstyle{J}},
!{(-1,-1.5)}*{-}
}
+(-1)^{\sigma_{P,Q}+|P||Q|+|P|}
\xygraph{
!{<0pt,0pt>;<10pt,0pt>:<0pt,10pt>::},
!{(-1,-1.5)}-!{(-2,-0.5)},
!{(-1,-1.5)}-!{(-1.5,-0.5)},
!{(-1,-1.5)}-!{(-1,-0.5)},
!{(-1,-1.5)}-!{(-0.5,-0.5)},
!{(-1,-1.5)}-!{(2,1.5)},
!{(0,-0.5)}-!{(-1,0.5)},
!{(0,-0.5)}-!{(-0.5,0.5)},
!{(0,-0.5)}-!{(0,0.5)},
!{(0,-0.5)}-!{(0.5,0.5)},
!{(1,0.5)}-!{(0,1.5)},
!{(1,0.5)}-!{(0.5,1.5)},
!{(1,0.5)}-!{(1,1.5)},
!{(1,0.5)}-!{(1.5,1.5)},
!{(-1.5,0)}*{\scriptstyle{P}},
!{(-0.5,1)}*{\scriptstyle{Q}},
!{(1,2)}*{\scriptstyle{J}},
!{(-1,-1.5)}*{-}
} \\
&+
\xygraph{
!{<0pt,0pt>;<10pt,0pt>:<0pt,10pt>::},
!{(-0.5,-0.5)}-!{(-0.5,-1.5)},
!{(-0.5,-0.5)}-!{(-1.5,0.5)},
!{(-0.5,-0.5)}-!{(-1,0.5)},
!{(-0.5,-0.5)}-!{(-0.5,0.5)},
!{(-0.5,-0.5)}-!{(0,0.5)},
!{(-0.5,-0.5)}-!{(0.5,0.5)},
!{(0,0.5)}-!{(-1,1.5)},
!{(0,0.5)}-!{(-0.5,1.5)},
!{(0,0.5)}-!{(0,1.5)},
!{(0,0.5)}-!{(0.5,1.5)},
!{(0,0.5)}-!{(1,1.5)},
!{(-0.5,-1.5)}-@(d,ur)!{(0.5,0.5)},
!{(-1,0.9)}*{\scriptstyle{I}},
!{(0,1.9)}*{\scriptstyle{J}}
}+(-1)^{|I|+1}\sum_{P\sqcup Q=J}(-1)^{\sigma_{P,Q}+|P||Q|}
\xygraph{
!{<0pt,0pt>;<10pt,0pt>:<0pt,10pt>::},
!{(-1,-1.5)}-!{(-2,-0.5)},
!{(-1,-1.5)}-!{(-1.5,-0.5)},
!{(-1,-1.5)}-!{(-1,-0.5)},
!{(-1,-1.5)}-!{(-0.5,-0.5)},
!{(-1,-1.5)}-!{(2,1.5)},
!{(0,-0.5)}-!{(-1,0.5)},
!{(0,-0.5)}-!{(-0.5,0.5)},
!{(0,-0.5)}-!{(0,0.5)},
!{(0,-0.5)}-!{(0.5,0.5)},
!{(1,0.5)}-!{(0,1.5)},
!{(1,0.5)}-!{(0.5,1.5)},
!{(1,0.5)}-!{(1,1.5)},
!{(1,0.5)}-!{(1.5,1.5)},
!{(-1.5,0)}*{\scriptstyle{I}},
!{(-0.5,1)}*{\scriptstyle{P}},
!{(1,2)}*{\scriptstyle{Q}},
!{(-1,-1.5)}*{-}
}.
\end{align*}

With the same reasoning, considering the sum
$$
\sum_{I'\sqcup J'=[n+1]}(-1)^{\sigma_{I',J'}+|I'||J'|}
\xygraph{
!{<0pt,0pt>;<10pt,0pt>:<0pt,10pt>::},
!{(-0.5,-0.5)}-!{(-0.5,-1.5)},
!{(-0.5,-0.5)}-!{(-1.5,0.5)},
!{(-0.5,-0.5)}-!{(-1,0.5)},
!{(-0.5,-0.5)}-!{(-0.5,0.5)},
!{(-0.5,-0.5)}-!{(0,0.5)},
!{(-0.5,-0.5)}-!{(1.5,1.5)},
!{(0.5,0.5)}-!{(-0.5,1.5)},
!{(0.5,0.5)}-!{(0,1.5)},
!{(0.5,0.5)}-!{(0.5,1.5)},
!{(0.5,0.5)}-!{(1,1.5)},
!{(-1,0.9)}*{\scriptstyle{I'}},
!{(0.5,1.9)}*{\scriptstyle{J'}}
},
$$
we get that
$$
\partial
\xygraph{
!{<0pt,0pt>;<10pt,0pt>:<0pt,10pt>::},
!{(0,0)}-!{(0,-1)},
!{(0,0)}-!{(-1,1)},
!{(0,0)}-!{(-0.5,1)},
!{(0,0)}-!{(0,1)},
!{(0,0)}-!{(0.5,1)},
!{(0,0)}-!{(1,1)},
!{(0,-1)}-@(d,ur)!{(1,1)},
!{(-1,1.4)}*{\scriptstyle{1}},
!{(0.5,1.4)}*{\scriptstyle{n}}
}=\sum_{I\sqcup J=[n]}(-1)^{\sigma_{I,J}+|I||J|+1}
\xygraph{
!{<0pt,0pt>;<10pt,0pt>:<0pt,10pt>::},
!{(-0.5,-0.5)}-!{(-0.5,-1.5)},
!{(-0.5,-0.5)}-!{(-1.5,0.5)},
!{(-0.5,-0.5)}-!{(-1,0.5)},
!{(-0.5,-0.5)}-!{(-0.5,0.5)},
!{(-0.5,-0.5)}-!{(0,0.5)},
!{(-0.5,-0.5)}-!{(0.5,0.5)},
!{(0,0.5)}-!{(-1,1.5)},
!{(0,0.5)}-!{(-0.5,1.5)},
!{(0,0.5)}-!{(0,1.5)},
!{(0,0.5)}-!{(0.5,1.5)},
!{(0,0.5)}-!{(1,1.5)},
!{(-0.5,-1.5)}-@(d,ur)!{(0.5,0.5)},
!{(-1,0.9)}*{\scriptstyle{I}},
!{(0,1.9)}*{\scriptstyle{J}}
}+(-1)^{\sigma_{I,J}+|I||J|+|I|}
\xygraph{
!{<0pt,0pt>;<10pt,0pt>:<0pt,10pt>::},
!{(-0.5,-0.5)}-!{(-0.5,-1.5)},
!{(-0.5,-0.5)}-!{(-1.5,0.5)},
!{(-0.5,-0.5)}-!{(-1,0.5)},
!{(-0.5,-0.5)}-!{(-0.5,0.5)},
!{(-0.5,-0.5)}-!{(0,0.5)},
!{(-0.5,-0.5)}-!{(1.5,1.5)},
!{(0.5,0.5)}-!{(-0.5,1.5)},
!{(0.5,0.5)}-!{(0,1.5)},
!{(0.5,0.5)}-!{(0.5,1.5)},
!{(0.5,0.5)}-!{(1,1.5)},
!{(-0.5,-1.5)}-@(d,ur)!{(1.5,1.5)},
!{(-1,0.9)}*{\scriptstyle{I}},
!{(0.5,1.9)}*{\scriptstyle{J}},
}.
$$
However, considering $I'=J$ and $J'=I$, since
\begin{align*}
(-1)^{\sigma_{I',J'}}=&(-1)^{\sigma_{J,I}}=(-1)^{|I||J|+\sigma_{I,J}},
\end{align*}
we have
$$
(-1)^{\sigma_{I',J'}+|I'||J'|+|I'|}
\xygraph{
!{<0pt,0pt>;<10pt,0pt>:<0pt,10pt>::},
!{(-0.5,-0.5)}-!{(-0.5,-1.5)},
!{(-0.5,-0.5)}-!{(-1.5,0.5)},
!{(-0.5,-0.5)}-!{(-1,0.5)},
!{(-0.5,-0.5)}-!{(-0.5,0.5)},
!{(-0.5,-0.5)}-!{(0,0.5)},
!{(-0.5,-0.5)}-!{(1.5,1.5)},
!{(0.5,0.5)}-!{(-0.5,1.5)},
!{(0.5,0.5)}-!{(0,1.5)},
!{(0.5,0.5)}-!{(0.5,1.5)},
!{(0.5,0.5)}-!{(1,1.5)},
!{(-0.5,-1.5)}-@(d,ur)!{(1.5,1.5)},
!{(-1,0.9)}*{\scriptstyle{I'}},
!{(0.5,1.9)}*{\scriptstyle{J'}},
}=
(-1)^{\sigma_{I,J}+|J|}
\xygraph{
!{<0pt,0pt>;<10pt,0pt>:<0pt,10pt>::},
!{(-0.5,-0.5)}-!{(-0.5,-1.5)},
!{(-0.5,-0.5)}-!{(-1.5,0.5)},
!{(-0.5,-0.5)}-!{(-1,0.5)},
!{(-0.5,-0.5)}-!{(-0.5,0.5)},
!{(-0.5,-0.5)}-!{(0,0.5)},
!{(-0.5,-0.5)}-!{(1.5,1.5)},
!{(0.5,0.5)}-!{(-0.5,1.5)},
!{(0.5,0.5)}-!{(0,1.5)},
!{(0.5,0.5)}-!{(0.5,1.5)},
!{(0.5,0.5)}-!{(1,1.5)},
!{(-0.5,-1.5)}-@(d,ur)!{(1.5,1.5)},
!{(-1,0.9)}*{\scriptstyle{J}},
!{(0.5,1.9)}*{\scriptstyle{I}},
}=
(-1)^{\sigma_{I,J}+|I||J|+|I|+1}
\xygraph{
!{<0pt,0pt>;<10pt,0pt>:<0pt,10pt>::},
!{(-0.5,-0.5)}-!{(-0.5,-1.5)},
!{(-0.5,-0.5)}-!{(-1.5,0.5)},
!{(-0.5,-0.5)}-!{(-1,0.5)},
!{(-0.5,-0.5)}-!{(-0.5,0.5)},
!{(-0.5,-0.5)}-!{(0,0.5)},
!{(-0.5,-0.5)}-!{(1.5,1.5)},
!{(0.5,0.5)}-!{(-0.5,1.5)},
!{(0.5,0.5)}-!{(0,1.5)},
!{(0.5,0.5)}-!{(0.5,1.5)},
!{(0.5,0.5)}-!{(1,1.5)},
!{(-0.5,-1.5)}-@(d,ur)!{(1.5,1.5)},
!{(-1,0.9)}*{\scriptstyle{I}},
!{(0.5,1.9)}*{\scriptstyle{J}},
},
$$
and the terms cancel pair-wise so that
$$
\partial
\xygraph{
!{<0pt,0pt>;<10pt,0pt>:<0pt,10pt>::},
!{(0,0)}-!{(0,-1)},
!{(0,0)}-!{(-1,1)},
!{(0,0)}-!{(-0.5,1)},
!{(0,0)}-!{(0,1)},
!{(0,0)}-!{(0.5,1)},
!{(0,0)}-!{(1,1)},
!{(0,-1)}-@(d,ur)!{(1,1)},
!{(-1,1.4)}*{\scriptstyle{1}},
!{(0.5,1.4)}*{\scriptstyle{n}}
}=\sum_{I\sqcup J=[n]}(-1)^{\sigma_{I,J}+|I||J|+1}
\xygraph{
!{<0pt,0pt>;<10pt,0pt>:<0pt,10pt>::},
!{(-0.5,-0.5)}-!{(-0.5,-1.5)},
!{(-0.5,-0.5)}-!{(-1.5,0.5)},
!{(-0.5,-0.5)}-!{(-1,0.5)},
!{(-0.5,-0.5)}-!{(-0.5,0.5)},
!{(-0.5,-0.5)}-!{(0,0.5)},
!{(-0.5,-0.5)}-!{(0.5,0.5)},
!{(0,0.5)}-!{(-1,1.5)},
!{(0,0.5)}-!{(-0.5,1.5)},
!{(0,0.5)}-!{(0,1.5)},
!{(0,0.5)}-!{(0.5,1.5)},
!{(0,0.5)}-!{(1,1.5)},
!{(-0.5,-1.5)}-@(d,ur)!{(0.5,0.5)},
!{(-1,0.9)}*{\scriptstyle{I}},
!{(0,1.9)}*{\scriptstyle{J}}
}.
$$

What remains to show is the vanishing of
\begin{align*}
\sum_{I\sqcup J=[n]}&\sum_{P\sqcup Q=I}(-1)^{\sigma_{I,J}+\sigma_{P,Q}+|I||J|+|P||Q|+1}
\xygraph{
!{<0pt,0pt>;<10pt,0pt>:<0pt,10pt>::},
!{(-1,-1.5)}-!{(-2,-0.5)},
!{(-1,-1.5)}-!{(-1.5,-0.5)},
!{(-1,-1.5)}-!{(-1,-0.5)},
!{(-1,-1.5)}-!{(0,0.5)},
!{(-1,-1.5)}-!{(2,1.5)},
!{(-0.5,-0.5)}-!{(-1.5,0.5)},
!{(-0.5,-0.5)}-!{(-1,0.5)},
!{(-0.5,-0.5)}-!{(-0.5,0.5)},
!{(-0.5,-0.5)}-!{(0.5,0.5)},
!{(1,0.5)}-!{(0,1.5)},
!{(1,0.5)}-!{(0.5,1.5)},
!{(1,0.5)}-!{(1,1.5)},
!{(1,0.5)}-!{(1.5,1.5)},
!{(-2,0)}*{\scriptstyle{P}},
!{(-0.5,1)}*{\scriptstyle{Q}},
!{(1,2)}*{\scriptstyle{J}},
!{(-1,-1.5)}*{-}
} \\
&+\sum_{I\sqcup J=[n]}\sum_{P\sqcup Q=I}(-1)^{\sigma_{I,J}+\sigma_{P,Q}+|I||J|+|P||Q|+|P|}
\xygraph{
!{<0pt,0pt>;<10pt,0pt>:<0pt,10pt>::},
!{(-1,-1.5)}-!{(-2,-0.5)},
!{(-1,-1.5)}-!{(-1.5,-0.5)},
!{(-1,-1.5)}-!{(-1,-0.5)},
!{(-1,-1.5)}-!{(-0.5,-0.5)},
!{(-1,-1.5)}-!{(2,1.5)},
!{(0,-0.5)}-!{(-1,0.5)},
!{(0,-0.5)}-!{(-0.5,0.5)},
!{(0,-0.5)}-!{(0,0.5)},
!{(0,-0.5)}-!{(0.5,0.5)},
!{(1,0.5)}-!{(0,1.5)},
!{(1,0.5)}-!{(0.5,1.5)},
!{(1,0.5)}-!{(1,1.5)},
!{(1,0.5)}-!{(1.5,1.5)},
!{(-1.5,0)}*{\scriptstyle{P}},
!{(-0.5,1)}*{\scriptstyle{Q}},
!{(1,2)}*{\scriptstyle{J}},
!{(-1,-1.5)}*{-}
} \\
&+\sum_{I\sqcup J=[n]}\sum_{P\sqcup Q=J}(-1)^{\sigma_{I,J}+\sigma_{P,Q}+|I||J|+|I|+|P||Q|+1}
\xygraph{
!{<0pt,0pt>;<10pt,0pt>:<0pt,10pt>::},
!{(-1,-1.5)}-!{(-2,-0.5)},
!{(-1,-1.5)}-!{(-1.5,-0.5)},
!{(-1,-1.5)}-!{(-1,-0.5)},
!{(-1,-1.5)}-!{(-0.5,-0.5)},
!{(-1,-1.5)}-!{(2,1.5)},
!{(0,-0.5)}-!{(-1,0.5)},
!{(0,-0.5)}-!{(-0.5,0.5)},
!{(0,-0.5)}-!{(0,0.5)},
!{(0,-0.5)}-!{(0.5,0.5)},
!{(1,0.5)}-!{(0,1.5)},
!{(1,0.5)}-!{(0.5,1.5)},
!{(1,0.5)}-!{(1,1.5)},
!{(1,0.5)}-!{(1.5,1.5)},
!{(-1.5,0)}*{\scriptstyle{I}},
!{(-0.5,1)}*{\scriptstyle{P}},
!{(1,2)}*{\scriptstyle{Q}},
!{(-1,-1.5)}*{-}
}.
\end{align*}

To see that the first sum in this expression equals zero, let $I'=P\sqcup J$, $J'=Q$, $P'=P$ and $Q'=J$, so that $I'=P\sqcup J$. We get
\begin{align*}
(-1)^{\sigma_{I',J'}+\sigma_{P',Q'}}=&(-1)^{\sigma_{I',Q}+\sigma_{P,J}}=(-1)^{\sigma_{P,J,Q}}=(-1)^{|Q||J|+\sigma_{P,Q,J}}=(-1)^{|Q||J|+\sigma_{I,J}+\sigma_{P,Q}}
\end{align*}
and hence
\begin{align*}
(-1)&^{\sigma_{I',J'}+\sigma_{P',Q'}+|I'||J'|+|P'||Q'|+1}
\xygraph{
!{<0pt,0pt>;<10pt,0pt>:<0pt,10pt>::},
!{(-1,-1.5)}-!{(-2,-0.5)},
!{(-1,-1.5)}-!{(-1.5,-0.5)},
!{(-1,-1.5)}-!{(-1,-0.5)},
!{(-1,-1.5)}-!{(0,0.5)},
!{(-1,-1.5)}-!{(2,1.5)},
!{(-0.5,-0.5)}-!{(-1.5,0.5)},
!{(-0.5,-0.5)}-!{(-1,0.5)},
!{(-0.5,-0.5)}-!{(-0.5,0.5)},
!{(-0.5,-0.5)}-!{(0.5,0.5)},
!{(1,0.5)}-!{(0,1.5)},
!{(1,0.5)}-!{(0.5,1.5)},
!{(1,0.5)}-!{(1,1.5)},
!{(1,0.5)}-!{(1.5,1.5)},
!{(-2,0)}*{\scriptstyle{P'}},
!{(-0.5,1)}*{\scriptstyle{Q'}},
!{(1,2)}*{\scriptstyle{J'}},
!{(-1,-1.5)}*{-}
}
=
(-1)^{\sigma_{I,J}+\sigma_{P,Q}+|Q||J|+|I'||Q|+|P||J|+1}
\xygraph{
!{<0pt,0pt>;<10pt,0pt>:<0pt,10pt>::},
!{(-1,-1.5)}-!{(-2,-0.5)},
!{(-1,-1.5)}-!{(-1.5,-0.5)},
!{(-1,-1.5)}-!{(-1,-0.5)},
!{(-1,-1.5)}-!{(0,0.5)},
!{(-1,-1.5)}-!{(2,1.5)},
!{(-0.5,-0.5)}-!{(-1.5,0.5)},
!{(-0.5,-0.5)}-!{(-1,0.5)},
!{(-0.5,-0.5)}-!{(-0.5,0.5)},
!{(-0.5,-0.5)}-!{(0.5,0.5)},
!{(1,0.5)}-!{(0,1.5)},
!{(1,0.5)}-!{(0.5,1.5)},
!{(1,0.5)}-!{(1,1.5)},
!{(1,0.5)}-!{(1.5,1.5)},
!{(-2,0)}*{\scriptstyle{P}},
!{(-0.5,1)}*{\scriptstyle{J}},
!{(1,2)}*{\scriptstyle{Q}},
!{(-1,-1.5)}*{-}
} \\
=&(-1)^{|Q||J|+\sigma_{I,J}+\sigma_{P,Q}+|I'||Q|+|P||J|}
\xygraph{
!{<0pt,0pt>;<10pt,0pt>:<0pt,10pt>::},
!{(-1,-1.5)}-!{(-2,-0.5)},
!{(-1,-1.5)}-!{(-1.5,-0.5)},
!{(-1,-1.5)}-!{(-1,-0.5)},
!{(-1,-1.5)}-!{(0.5,1.5)},
!{(-1,-1.5)}-!{(1.5,-0.5)},
!{(0,0.5)}-!{(-1,1.5)},
!{(0,0.5)}-!{(-0.5,1.5)},
!{(0,0.5)}-!{(0,1.5)},
!{(0,0.5)}-!{(1,1.5)},
!{(1.5,-0.5)}-!{(0.5,0.5)},
!{(1.5,-0.5)}-!{(1,0.5)},
!{(1.5,-0.5)}-!{(1.5,0.5)},
!{(1.5,-0.5)}-!{(2,0.5)},
!{(1.5,-0.5)}-!{(2.5,0.5)},
!{(-1.5,0)}*{\scriptstyle{P}},
!{(1.5,1)}*{\scriptstyle{J}},
!{(0,2)}*{\scriptstyle{Q}},
!{(-1,-1.5)}*{-}
}
=(-1)^{\sigma_{I,J}+\sigma_{P,Q}+|I'||Q|+|P||J|}
\xygraph{
!{<0pt,0pt>;<10pt,0pt>:<0pt,10pt>::},
!{(-1,-1.5)}-!{(-2,-0.5)},
!{(-1,-1.5)}-!{(-1.5,-0.5)},
!{(-1,-1.5)}-!{(-1,-0.5)},
!{(-1,-1.5)}-!{(0,0.5)},
!{(-1,-1.5)}-!{(2,1.5)},
!{(-0.5,-0.5)}-!{(-1.5,0.5)},
!{(-0.5,-0.5)}-!{(-1,0.5)},
!{(-0.5,-0.5)}-!{(-0.5,0.5)},
!{(-0.5,-0.5)}-!{(0.5,0.5)},
!{(1,0.5)}-!{(0,1.5)},
!{(1,0.5)}-!{(0.5,1.5)},
!{(1,0.5)}-!{(1,1.5)},
!{(1,0.5)}-!{(1.5,1.5)},
!{(-2,0)}*{\scriptstyle{P}},
!{(-0.5,1)}*{\scriptstyle{Q}},
!{(1,2)}*{\scriptstyle{J}},
!{(-1,-1.5)}*{-}
} \\
=&(-1)^{\sigma_{I,J}+\sigma_{P,Q}+|P||Q|+|J||Q|+|P||J|}
\xygraph{
!{<0pt,0pt>;<10pt,0pt>:<0pt,10pt>::},
!{(-1,-1.5)}-!{(-2,-0.5)},
!{(-1,-1.5)}-!{(-1.5,-0.5)},
!{(-1,-1.5)}-!{(-1,-0.5)},
!{(-1,-1.5)}-!{(0,0.5)},
!{(-1,-1.5)}-!{(2,1.5)},
!{(-0.5,-0.5)}-!{(-1.5,0.5)},
!{(-0.5,-0.5)}-!{(-1,0.5)},
!{(-0.5,-0.5)}-!{(-0.5,0.5)},
!{(-0.5,-0.5)}-!{(0.5,0.5)},
!{(1,0.5)}-!{(0,1.5)},
!{(1,0.5)}-!{(0.5,1.5)},
!{(1,0.5)}-!{(1,1.5)},
!{(1,0.5)}-!{(1.5,1.5)},
!{(-2,0)}*{\scriptstyle{P}},
!{(-0.5,1)}*{\scriptstyle{Q}},
!{(1,2)}*{\scriptstyle{J}},
!{(-1,-1.5)}*{-}
}
=(-1)^{\sigma_{I,J}+\sigma_{P,Q}+|P||Q|+|I||J|}
\xygraph{
!{<0pt,0pt>;<10pt,0pt>:<0pt,10pt>::},
!{(-1,-1.5)}-!{(-2,-0.5)},
!{(-1,-1.5)}-!{(-1.5,-0.5)},
!{(-1,-1.5)}-!{(-1,-0.5)},
!{(-1,-1.5)}-!{(0,0.5)},
!{(-1,-1.5)}-!{(2,1.5)},
!{(-0.5,-0.5)}-!{(-1.5,0.5)},
!{(-0.5,-0.5)}-!{(-1,0.5)},
!{(-0.5,-0.5)}-!{(-0.5,0.5)},
!{(-0.5,-0.5)}-!{(0.5,0.5)},
!{(1,0.5)}-!{(0,1.5)},
!{(1,0.5)}-!{(0.5,1.5)},
!{(1,0.5)}-!{(1,1.5)},
!{(1,0.5)}-!{(1.5,1.5)},
!{(-2,0)}*{\scriptstyle{P}},
!{(-0.5,1)}*{\scriptstyle{Q}},
!{(1,2)}*{\scriptstyle{J}},
!{(-1,-1.5)}*{-}
},
\end{align*}
so that the terms cancel pair-wise. Next, in the fourth sum, consider $I'=P$, $J'=Q\sqcup J$, $P'=Q$ and $Q'=J$. We have
\begin{align*}
(-1)^{\sigma_{I',J'}+\sigma_{P',Q'}}=(-1)^{\sigma_{P,J'}+\sigma_{Q,J}}=(-1)^{\sigma_{P,Q,J}}=(-1)^{\sigma_{I,J}+\sigma_{P,Q}},
\end{align*}
so that
\begin{align*}
(-1)&^{\sigma_{I',J'}+\sigma_{P',Q'}+|I'||J'|+|I'|+|P'||Q'|+1}
\xygraph{
!{<0pt,0pt>;<10pt,0pt>:<0pt,10pt>::},
!{(-1,-1.5)}-!{(-2,-0.5)},
!{(-1,-1.5)}-!{(-1.5,-0.5)},
!{(-1,-1.5)}-!{(-1,-0.5)},
!{(-1,-1.5)}-!{(-0.5,-0.5)},
!{(-1,-1.5)}-!{(2,1.5)},
!{(0,-0.5)}-!{(-1,0.5)},
!{(0,-0.5)}-!{(-0.5,0.5)},
!{(0,-0.5)}-!{(0,0.5)},
!{(0,-0.5)}-!{(0.5,0.5)},
!{(1,0.5)}-!{(0,1.5)},
!{(1,0.5)}-!{(0.5,1.5)},
!{(1,0.5)}-!{(1,1.5)},
!{(1,0.5)}-!{(1.5,1.5)},
!{(-1.5,0)}*{\scriptstyle{I'}},
!{(-0.5,1)}*{\scriptstyle{P'}},
!{(1,2)}*{\scriptstyle{Q'}},
!{(-1,-1.5)}*{-}
}
=(-1)^{\sigma_{I,J}+\sigma_{P,Q}+|P|(|Q|+|J|)+|P|+|Q||J|+1}
\xygraph{
!{<0pt,0pt>;<10pt,0pt>:<0pt,10pt>::},
!{(-1,-1.5)}-!{(-2,-0.5)},
!{(-1,-1.5)}-!{(-1.5,-0.5)},
!{(-1,-1.5)}-!{(-1,-0.5)},
!{(-1,-1.5)}-!{(-0.5,-0.5)},
!{(-1,-1.5)}-!{(2,1.5)},
!{(0,-0.5)}-!{(-1,0.5)},
!{(0,-0.5)}-!{(-0.5,0.5)},
!{(0,-0.5)}-!{(0,0.5)},
!{(0,-0.5)}-!{(0.5,0.5)},
!{(1,0.5)}-!{(0,1.5)},
!{(1,0.5)}-!{(0.5,1.5)},
!{(1,0.5)}-!{(1,1.5)},
!{(1,0.5)}-!{(1.5,1.5)},
!{(-1.5,0)}*{\scriptstyle{P}},
!{(-0.5,1)}*{\scriptstyle{Q}},
!{(1,2)}*{\scriptstyle{J}},
!{(-1,-1.5)}*{-}
} \\
=&(-1)^{\sigma_{I,J}+\sigma_{P,Q}+|P||Q|+|I||J|+|P|+1}
\xygraph{
!{<0pt,0pt>;<10pt,0pt>:<0pt,10pt>::},
!{(-1,-1.5)}-!{(-2,-0.5)},
!{(-1,-1.5)}-!{(-1.5,-0.5)},
!{(-1,-1.5)}-!{(-1,-0.5)},
!{(-1,-1.5)}-!{(-0.5,-0.5)},
!{(-1,-1.5)}-!{(2,1.5)},
!{(0,-0.5)}-!{(-1,0.5)},
!{(0,-0.5)}-!{(-0.5,0.5)},
!{(0,-0.5)}-!{(0,0.5)},
!{(0,-0.5)}-!{(0.5,0.5)},
!{(1,0.5)}-!{(0,1.5)},
!{(1,0.5)}-!{(0.5,1.5)},
!{(1,0.5)}-!{(1,1.5)},
!{(1,0.5)}-!{(1.5,1.5)},
!{(-1.5,0)}*{\scriptstyle{P}},
!{(-0.5,1)}*{\scriptstyle{Q}},
!{(1,2)}*{\scriptstyle{J}},
!{(-1,-1.5)}*{-}
}.
\end{align*}
Hence, the fourth sum cancels the second.
\epro

\bth
The wheeled operad $\ULie$ is wheeled Koszul.
\eth
\bpro
We know that $H(\Omega^\cal(\ULie^\text{!`})_o)=H(\Lie_\infty)=\Lie$. Hence, we need to show that the cohomology of $C=\Omega^\circlearrowleft(\ULie^\text{!`})_w$ equals $\ULie_w$.

The graphs of genus one span a subcomplex of $C$ isomorphic to $((\Lie_\infty)^\cal)_w$ which, since $\Lie$ is stably Koszul, equals $(\Lie_\infty^\cal)_w$. Considering the inclusion of $(\Lie_\infty^\cal)_w$ into $C$, we see that if we show that the cohomology of $C_0=C/(\Lie_\infty^\cal)_w$ equals
$\k
\xygraph{
!{<0pt,0pt>;<10pt,0pt>:<0pt,10pt>::},
!{(0,0)}-!{(0,0.5)},
!{(0,0)}*{-}
}$,
then we are done, because then the long exact sequence in cohomology reduces to
$$
0 \ra\k
\xygraph{
!{<0pt,0pt>;<10pt,0pt>:<0pt,10pt>::},
!{(0,0)}-!{(0,0.5)},
!{(0,0)}*{-}
}
\ra\Lie_w^\cal\ra H^0(C)\ra 0
$$
and
$$
H(C)=\frac{\Lie^\circlearrowleft_w}{
\xygraph{
!{<0pt,0pt>;<5pt,0pt>:<0pt,5pt>::},
!{(0,0)}-!{(0,-1)},
!{(0,0)}-!{(-1,1)},
!{(0,0)}-!{(1,1)},
!{(0,-1)}-@(d,ur)!{(1,1)}
}}=\ULie_w.
$$
Now, the subcomplex of $C_0$ consisting of decorated trees without root leg whose root vertex has only one incoming leg is isomorphic to $\Lie_\infty[1]$ and the quotient is isomorphic to $\overline{\Lie_\infty}[2]$. Considering the long exact sequence in cohomology induced by the inclusion of this subcomplex into $C_0$, that is
$$
0\ra H^{-2}(C_0)\ra\overline{\Lie}[2]\ra\k
\xygraph{
!{<0pt,0pt>;<10pt,0pt>:<0pt,10pt>::},
!{(0,0)}-!{(0,0.5)},
!{(0,0)}*{-}
}
\oplus\overline{\Lie}[1]\ra H^{-1}(C_0)\ra 0,
$$
we conclude that
$H(C_0)=\k
\xygraph{
!{<0pt,0pt>;<10pt,0pt>:<0pt,10pt>::},
!{(0,0)}-!{(0,0.5)},
!{(0,0)}*{-}
}$.
\epro

\bco
The minimal resolution $\ULie_\infty$ equals $\Omega^\cal(\ULie^\text{!`})$ as given in Proposition \ref{prop:threeone}.
\eco

\bre
The proof that
$H(C_0)=\k
\xygraph{
!{<0pt,0pt>;<10pt,0pt>:<0pt,10pt>::},
!{(0,0)}-!{(0,0.5)},
!{(0,0)}*{-}
}$
is essentially the same as that of the claim in the proof of Theorem $4.1.1$ in \cite{Merkulov2007}. As a second remark we note that since $\ULie^!=\Com^\circlearrowleft$, the above proof gives a new and shorter proof of the theorem that $\Com^\circlearrowleft$ is wheeled Koszul proven in \cite{Markl2007}.
\ere

\bde
A structure of \emph{unimodular $L_\infty$-algebra} on a finite dimensional dg vector space $V$ is a representation of $\ULie_\infty$ in $V$. By the corollary above, this is equivalent with a finite dimensional dg vector space $V$ together with two families of morphisms
$$
\{\ell_n:\wedge^n V\ra V\}_{n=2}^\infty
$$
with $|\ell_n|=2-n$ and
$$
\{q_n:\wedge^n V\ra\k\}_{n=1}^\infty
$$
with $|q_n|=-n$, satisfying the equations
$$
\partial\ell_n=\sum_{I\sqcup J=[n]}(-1)^{\sigma_{I,J}+|I||J|}\ell_{|I|+1}(\id\otimes\cdots\otimes\id\otimes\ell_{|J|})\sigma_{I,J}
$$
and
$$
\partial q_n=\sum_{I\sqcup J=[n]}(-1)^{\sigma_{I,J}+|I||J|}q_{|I|+1}(\id\otimes\cdots\otimes\id\otimes\ell_{|J|})\sigma_{I,J}+\tr_{n+1}\ell_{n+1}.
$$
\ede

\section{Transfer of structure}

Assume we have a diagram
$$
\xymatrix@1{
V \ar@/^/[r]^{g} \ar@(ul,dl)_{h} & U \ar@/^/[l]^{f}
}
$$
where $V$ is a unimodular $L_\infty$-algebra, $U$ is a finite dimensional dg vector space and $h$ is a cochain homotopy between $fg$ and $\id_V$, i.e.~a degree $-1$ map satisfying
$$
fg-\id_V=d_V h+hd_V.
$$
In this situation we give, explicitly, an induced structure of unimodular $L_\infty$-algebra on $U$.

To handle signs we apply the parity change functor to $\ULie_\infty$, $\Pi\ULie_\infty$ is the free wheeled operad on symmetric generators
$$
\left\{
\xygraph{
!{<0pt,0pt>;<10pt,0pt>:<0pt,10pt>::},
!{(0,0)}-!{(0,-1)},
!{(0,0)}-!{(-1,1)},
!{(0,0)}-!{(-0.5,1)},
!{(0,0)}-!{(0,1)},
!{(0,0)}-!{(0.5,1)},
!{(0,0)}-!{(1,1)},
!{(-1,1.4)}*{\scriptstyle{1}},
!{(1,1.4)}*{\scriptstyle{n}}
}
\right\}_{n=2}^\infty
\text{ and }
\left\{
\xygraph{
!{<0pt,0pt>;<10pt,0pt>:<0pt,10pt>::},
!{(0,0)}-!{(-1,1)},
!{(0,0)}-!{(-0.5,1)},
!{(0,0)}-!{(0,1)},
!{(0,0)}-!{(0.5,1)},
!{(0,0)}-!{(1,1)},
!{(-1,1.4)}*{\scriptstyle{1}},
!{(1,1.4)}*{\scriptstyle{n}},
!{(0,0)}*{-}
}
\right\}_{n=1}^\infty,
$$
of degrees $1$ and $0$ respectively. The differential is given on generators by
\begin{align}\label{eq:fourone}
\partial
\xygraph{
!{<0pt,0pt>;<10pt,0pt>:<0pt,10pt>::},
!{(0,0)}-!{(0,-1)},
!{(0,0)}-!{(-1,1)},
!{(0,0)}-!{(-0.5,1)},
!{(0,0)}-!{(0,1)},
!{(0,0)}-!{(0.5,1)},
!{(0,0)}-!{(1,1)},
!{(-1,1.4)}*{\scriptstyle{1}},
!{(1,1.4)}*{\scriptstyle{n}}
}=&
-\sum_{\substack{I\sqcup J=[n] \\ 1\leq |I|\leq n-2}}
\xygraph{
!{<0pt,0pt>;<10pt,0pt>:<0pt,10pt>::},
!{(-0.5,-0.5)}-!{(-0.5,-1.5)},
!{(-0.5,-0.5)}-!{(-1.5,0.5)},
!{(-0.5,-0.5)}-!{(-1,0.5)},
!{(-0.5,-0.5)}-!{(-0.5,0.5)},
!{(-0.5,-0.5)}-!{(0,0.5)},
!{(-0.5,-0.5)}-!{(1.5,1.5)},
!{(0.5,0.5)}-!{(-0.5,1.5)},
!{(0.5,0.5)}-!{(0,1.5)},
!{(0.5,0.5)}-!{(0.5,1.5)},
!{(0.5,0.5)}-!{(1,1.5)},
!{(-1,0.9)}*{\scriptstyle{I}},
!{(0.5,1.9)}*{\scriptstyle{J}}
}, \\
\label{eq:fourtwo}
\partial
\xygraph{
!{<0pt,0pt>;<10pt,0pt>:<0pt,10pt>::},
!{(0,0)}-!{(-1,1)},
!{(0,0)}-!{(-0.5,1)},
!{(0,0)}-!{(0,1)},
!{(0,0)}-!{(0.5,1)},
!{(0,0)}-!{(1,1)},
!{(-1,1.4)}*{\scriptstyle{1}},
!{(1,1.4)}*{\scriptstyle{n}},
!{(0,0)}*{-}
}=&
\sum_{\substack{I\sqcup J=[n] \\ 0\leq |I|\leq n-2}}
\xygraph{
!{<0pt,0pt>;<10pt,0pt>:<0pt,10pt>::},
!{(-0.5,-0.5)}-!{(-1.5,0.5)},
!{(-0.5,-0.5)}-!{(-1,0.5)},
!{(-0.5,-0.5)}-!{(-0.5,0.5)},
!{(-0.5,-0.5)}-!{(0,0.5)},
!{(-0.5,-0.5)}-!{(1.5,1.5)},
!{(0.5,0.5)}-!{(-0.5,1.5)},
!{(0.5,0.5)}-!{(0,1.5)},
!{(0.5,0.5)}-!{(0.5,1.5)},
!{(0.5,0.5)}-!{(1,1.5)},
!{(-1,0.9)}*{\scriptstyle{I}},
!{(0.5,1.9)}*{\scriptstyle{J}},
!{(-0.5,-0.5)}*{-}
}+
\xygraph{
!{<0pt,0pt>;<10pt,0pt>:<0pt,10pt>::},
!{(0,0)}-!{(0,-1)},
!{(0,0)}-!{(-1,1)},
!{(0,0)}-!{(-0.5,1)},
!{(0,0)}-!{(0,1)},
!{(0,0)}-!{(0.5,1)},
!{(0,0)}-!{(1,1)},
!{(0,-1)}-@(d,ur)!{(1,1)},
!{(-1,1.4)}*{\scriptstyle{1}},
!{(0.5,1.4)}*{\scriptstyle{n}}
}.
\end{align}

Let $\fG_1$ denote the class of directed trees with root leg, $\fG_0$ the class of directed trees without root leg and $\fG_1^\cal$ the class of graphs obtained by grafting the root leg of a tree to one of its ingoing legs. Let $\fG_1(n)$ denote the subclass of trees with $n$ ingoing legs and similarly for $\fG_0(n)$ and $\fG_1^\cal(n)$. For any tree $G\in\fG_1(n)$ we define a map $\theta_G:V[1]^{\otimes n}\ra V[1]$ by decorating the root vertex $v$ with $\tilde{\ell}_{|in(v)|}$ and every other vertex $v$ with $h\tilde{\ell}_{|in(v)|}$. For any tree $G\in\fG_0(n)$ we define a morphism $\theta_G:V[1]^{\otimes n}\ra\k$ by decorating the root vertex $v$ with $\tilde{q}_{|in(v)|}$ and every other vertex $v$ with $h\tilde{\ell}_{|in(v)|}$. For a graph $G\in\fG^\cal_1(n)$, we define $\theta_G:V[1]^{\otimes n}\ra\k$ by decorating all vertices $v$ with $h\tilde{\ell}_{|in(v)|}$. Next we define a morphism $\Pi\ULie_\infty\ra\oEnd^\cal_{U[1]}$ by
\begin{align*}
\xygraph{
!{<0pt,0pt>;<10pt,0pt>:<0pt,10pt>::},
!{(0,0)}-!{(0,-1)},
!{(0,0)}-!{(-1,1)},
!{(0,0)}-!{(-0.5,1)},
!{(0,0)}-!{(0,1)},
!{(0,0)}-!{(0.5,1)},
!{(0,0)}-!{(1,1)},
!{(-1,1.4)}*{\scriptstyle{1}},
!{(1,1.4)}*{\scriptstyle{n}}
}
\mapsto&\sum_{G\in\fG_1(n)}g\theta_G f^{\otimes n}, \\
\xygraph{
!{<0pt,0pt>;<10pt,0pt>:<0pt,10pt>::},
!{(0,0)}-!{(-1,1)},
!{(0,0)}-!{(-0.5,1)},
!{(0,0)}-!{(0,1)},
!{(0,0)}-!{(0.5,1)},
!{(0,0)}-!{(1,1)},
!{(-1,1.4)}*{\scriptstyle{1}},
!{(1,1.4)}*{\scriptstyle{n}},
!{(0,0)}*{-}
}
\mapsto&\sum_{G\in(\fG_0\sqcup\fG_1^\cal)(n)}\theta_G f^{\otimes n}.
\end{align*}

\bth
The above morphism $\Pi\ULie_\infty\ra\oEnd^\cal_{U[1]}$ is a morphism of dg wheeled operads and hence defines a unimodular $L_\infty$-algebra structure on $U$.
\eth
\bpro
For a graph $G$ and an edge $e\in\mathbf{e}(G)$, let $\theta_{G,e}^\circ$ denote the morphism that differs from $\theta_G$ in the respect that the vertex $v$ from which $e$ starts is decorated with $fg\tilde{\ell}_{|in(v)|}$ instead of $h\tilde{\ell}_{|in(v)|}$. For the equality \eqref{eq:fourone} we need to show that
$$
\sum_{G\in\fG_1(n)}\partial g\theta_Gf^{\otimes n}=\sum_{G\in\fG_1(n)}dg\theta_G f^{\otimes n}+g\theta_G f^{\otimes n}d_\otimes
$$
equals
$$
-\sum_{\substack{G_1\in\fG_1(i+1) \\ G_2\in\fG_1(j) \\ n=i+j}}g\theta_{G_1}(f\otimes\cdots\otimes f\otimes fg\theta_{G_2}(f\otimes\cdots\otimes f))
=-\sum_{G\in\fG_1(n)}\sum_{e\in\mathbf{e}(G)}g\theta_{G,e}^\circ f^{\otimes n}.
$$

Since $f$ and $g$ commute with differentials, this would follow, by considering one tree at a time, from
$$
\theta_G d_\otimes=-d\theta_G+\sum_{e\in\mathbf{e}(G)}(-\theta_{G,e}^\circ+\theta_{G,e}^{\id})-\sum_{\substack{(G',e') \\ G'\in\fG_1(n) \\ G'/e'=G}}\theta_{G',e'}^{\id}.
$$
This in turn follows by commuting the differential in $\theta_G d_\otimes$ from the right to the left. Considering one vertex of $G$ at a time, we note that for a non-root vertex
$$
h\tilde{\ell_k}d_\otimes=-hd\tilde{\ell_k}-\sum h\tilde{\ell}_{i+1}(\id\otimes\cdots\otimes\tilde{\ell}_j)\sigma
=dh\tilde{\ell}_k-fg\tilde{\ell}_k+\tilde{\ell}_k-\sum h\tilde{\ell}_{i+1}(\id\otimes\cdots\otimes\tilde{\ell}_j)\sigma,
$$
while for the root vertex
$$
\tilde{\ell_k}d_\otimes=-d\tilde{\ell_k}-\sum\tilde{\ell}_{i+1}(\id\otimes\cdots\otimes\tilde{\ell}_j)\sigma.
$$

For the equality \eqref{eq:fourtwo} we need to show that
$$
\sum_{G\in(\fG_0\sqcup\fG_1^\cal)(n)}\partial\theta_G f^{\otimes n}=-\sum_{G\in(\fG_0\sqcup\fG_1^\cal)(n)}\theta_G f^{\otimes n}d_\otimes
$$
equals
\begin{align*}
\sum_{\substack{G_1\in(\fG_0\sqcup\fG_1^\cal)(i+1) \\ G_2\in\fG_1(j) \\ n=i+j}}&\theta_{G_1}(f\otimes\cdots\otimes f\otimes fg\theta_{G_2}(f\otimes\cdots\otimes f))
+\sum_{G\in\fG_1(n+1)}\tr_{n+1}(g\theta_G(f\otimes\cdots\otimes f)) \\
=&\left(\sum_{G\in\fG_0(n)}\sum_{e\in\mathbf{e}(G)}+\sum_{G\in\fG_1^\cal(n)}\sum_{e\in\mathbf{e_{nc}}(G)}+\sum_{G\in\fG_1^\cal(n)}\sum_{e\in\mathbf{e}_c(G)}\right)\theta_{G,e}^\circ f^{\otimes n} \\
=&\sum_{G\in(\fG_0\sqcup\fG_1^\cal)(n)}\sum_{e\in\mathbf{e}(G)}\theta_{G,e}^\circ f^{\otimes n}.
\end{align*}
This would follow from
$$
-\theta_G d_\otimes=\sum_{e\in\mathbf{e}(G)}(\theta_{G,e}^\circ-\theta_{G,e}^{\id})+\sum_{\substack{(G',e') \\ G'\in(\fG_0\sqcup\fG_1^\cal)(n) \\ G'/e'=G}}\theta_{G',e'}^{\id}
$$
which again follows from commuting the differential from right to left. Considering one vertex at a time, we have for a non-root, non-cyclic vertex that
$$
-h\tilde{\ell}_k d_\otimes=-dh\tilde{\ell}_k+fg\tilde{\ell}_k-\tilde{\ell}_k+\sum h\tilde{\ell}_{i+1}(\id\otimes\cdots\otimes\tilde{\ell}_j)\sigma.
$$
For the root vertex of a $G\in\fG_0(n)$ we have
$$
-\tilde{q}_k d_\otimes=\sum\tilde{q}_{i+1}(\id\otimes\cdots\otimes\tilde{\ell}_j)\sigma+\tr_{k+1}(\tilde{\ell}_{k+1}).
$$
Finally, for a cyclic vertex $v_k$ of a $G\in\fG_1^\cal(n)$,
\begin{align*}
-h\tilde{\ell}_k&((d\otimes\id\otimes\cdots\otimes\id)+\cdots+(\id\otimes\cdots\otimes d\otimes\id)) \\
=&h\tilde{\ell}_k(\id\otimes\cdots\otimes\id\otimes d)+hd\tilde{\ell}_k+\sum h\tilde{\ell}_{i+1}(\id\otimes\cdots\otimes\tilde{\ell}_j)\sigma.
\end{align*}
The first term gives a term $dh\tilde{\ell}_{k'}$ for the previous vertex $v_{k'}$ in the cycle and like-wise, the above equality for the next vertex in the cycle gives a term $dh\tilde{\ell}_k$ to our expression for $v_k$. Hence, we get
\begin{align*}
hd\tilde{\ell}_k+&\sum h\tilde{\ell}_{i+1}(\id\otimes\cdots\otimes\tilde{\ell}_j)\sigma+dh\tilde{\ell}_k \\
=&fg\tilde{\ell}_k-\tilde{\ell}_k+\sum h\tilde{\ell}_{i+1}(\id\otimes\cdots\otimes\tilde{\ell}_j)\sigma+dh\tilde{\ell}_k.
\end{align*}
\epro

\bex
As an illustration of the proof given above, we will show the required equality \eqref{eq:fourtwo} in the case $n=2$ for the induced structure graphically. The following is a legend for our graphical notations:
\begin{align*}
\tilde{q}_k=&
\xygraph{
!{<0pt,0pt>;<10pt,0pt>:<0pt,10pt>::},
!{(0,0)}-!{(-1,1)},
!{(0,0)}-!{(-0.5,1)},
!{(0,0)}-!{(0,1)},
!{(0,0)}-!{(0.5,1)},
!{(0,0)}-!{(1,1)},
!{(0,0)}*{-},
!{(-1,1.4)}*{\scriptstyle{1}},
!{(1,1.4)}*{\scriptstyle{k}}
}
&
\tilde{\ell}_k=&
\xygraph{
!{<0pt,0pt>;<10pt,0pt>:<0pt,10pt>::},
!{(0,0)}-!{(0,-1)},
!{(0,0)}-!{(-1,1)},
!{(0,0)}-!{(-0.5,1)},
!{(0,0)}-!{(0,1)},
!{(0,0)}-!{(0.5,1)},
!{(0,0)}-!{(1,1)},
!{(-1,1.4)}*{\scriptstyle{1}},
!{(1,1.4)}*{\scriptstyle{k}}
} \\
h=&
\xygraph{
!{<0pt,0pt>;<10pt,0pt>:<0pt,10pt>::},
!{(0,-0.5)}-!{(0,0.5)},
!{(0,0)}*{\scriptstyle{\bullet}}
}
&
d=&
\xygraph{
!{<0pt,0pt>;<10pt,0pt>:<0pt,10pt>::},
!{(0,-0.5)}-!{(0,0.5)},
!{(0,0)}*{\scriptstyle{\times}}
} \\
fg=&
\xygraph{
!{<0pt,0pt>;<10pt,0pt>:<0pt,10pt>::},
!{(0,-0.5)}-!{(0,0.5)},
!{(0,0)}*{\scriptstyle{\circ}}
}
&&
\end{align*}

Ignoring the $f$:s decorating the leaves, the left hand side of \eqref{eq:fourtwo} for the induced structure is
\begin{align*}
\partial\left(
\xygraph{
!{<0pt,0pt>;<10pt,0pt>:<0pt,10pt>::},
!{(-1,1)}-!{(0,0)}-!{(1,1)},
!{(0,0)}*{-}
}
+
\xygraph{
!{<0pt,0pt>;<10pt,0pt>:<0pt,10pt>::},
!{(-1,1)}-!{(0,0)}-!{(1,1)},
!{(0,0)}-!{(0,-1)},
!{(0,-0.5)}*{\scriptstyle{\bullet}},
!{(0,-1)}*{-}
}
+
\xygraph{
!{<0pt,0pt>;<10pt,0pt>:<0pt,10pt>::},
!{(0.5,-0.5)}-!{(0.5,-1.5)},
!{(0.5,-0.5)}-!{(-1.5,1.5)},
!{(0.5,-0.5)}-!{(1.5,0.5)},
!{(-0.5,0.5)}-!{(0.5,1.5)},
!{(0.5,-1.5)}-@(d,ur)!{(1.5,0.5)},
!{(0.5,-1)}*{\scriptstyle{\bullet}},
!{(0,0)}*{\scriptstyle{\bullet}}
}
+
\xygraph{
!{<0pt,0pt>;<10pt,0pt>:<0pt,10pt>::},
!{(-0.5,-0.5)}-!{(-0.5,-1.5)},
!{(-0.5,-0.5)}-!{(-1.5,0.5)},
!{(-0.5,-0.5)}-!{(1.5,1.5)},
!{(0.5,0.5)}-!{(-0.5,1.5)},
!{(-0.5,-1.5)}-@(d,ur)!{(1.5,1.5)},
!{(-0.5,-1)}*{\scriptstyle{\bullet}},
!{(0,0)}*{\scriptstyle{\bullet}}
}
+
\xygraph{
!{<0pt,0pt>;<10pt,0pt>:<0pt,10pt>::},
!{(0,0)}-!{(0,-1)},
!{(0,0)}-!{(-1,1)},
!{(0,0)}-!{(0,1)},
!{(0,0)}-!{(1,1)},
!{(0,-1)}-@(d,ur)!{(1,1)},
!{(0,-0.5)}*{\scriptstyle{\bullet}}
}
\right),
\end{align*}
while the right hand side is
\begin{align*}
\xygraph{
!{<0pt,0pt>;<10pt,0pt>:<0pt,10pt>::},
!{(0,0)}-!{(0,-1)},
!{(0,0)}-!{(-1,1)},
!{(0,0)}-!{(1,1)},
!{(0,-1)}*{-},
!{(0,-0.5)}*{\scriptstyle{\circ}},
}
+
\xygraph{
!{<0pt,0pt>;<10pt,0pt>:<0pt,10pt>::},
!{(0.5,-0.5)}-!{(0.5,-1.5)},
!{(0.5,-0.5)}-!{(-1.5,1.5)},
!{(0.5,-0.5)}-!{(1.5,0.5)},
!{(-0.5,0.5)}-!{(0.5,1.5)},
!{(0.5,-1.5)}-@(d,ur)!{(1.5,0.5)},
!{(0.5,-1)}*{\scriptstyle{\bullet}},
!{(0,0)}*{\scriptstyle{\circ}}
}
+
\xygraph{
!{<0pt,0pt>;<10pt,0pt>:<0pt,10pt>::},
!{(0,0)}-!{(0,-1)},
!{(0,0)}-!{(-1,1)},
!{(0,0)}-!{(0,1)},
!{(0,0)}-!{(1,1)},
!{(0,-1)}-@(d,ur)!{(1,1)},
!{(0,-0.5)}*{\scriptstyle{\circ}}
}
+
\xygraph{
!{<0pt,0pt>;<10pt,0pt>:<0pt,10pt>::},
!{(-0.5,-0.5)}-!{(-0.5,-1.5)},
!{(-0.5,-0.5)}-!{(-1.5,0.5)},
!{(-0.5,-0.5)}-!{(1.5,1.5)},
!{(0.5,0.5)}-!{(-0.5,1.5)},
!{(-0.5,-1.5)}-@(d,ur)!{(1.5,1.5)},
!{(-0.5,-1)}*{\scriptstyle{\circ}},
!{(0,0)}*{\scriptstyle{\bullet}}
}
+
\xygraph{
!{<0pt,0pt>;<10pt,0pt>:<0pt,10pt>::},
!{(-0.5,-0.5)}-!{(-0.5,-1.5)},
!{(-0.5,-0.5)}-!{(-1.5,0.5)},
!{(-0.5,-0.5)}-!{(1.5,1.5)},
!{(0.5,0.5)}-!{(-0.5,1.5)},
!{(-0.5,-1.5)}-@(d,ur)!{(1.5,1.5)},
!{(-0.5,-1)}*{\scriptstyle{\bullet}},
!{(0,0)}*{\scriptstyle{\circ}}
}
+
\xygraph{
!{<0pt,0pt>;<10pt,0pt>:<0pt,10pt>::},
!{(0.5,-0.5)}-!{(0.5,-1.5)},
!{(0.5,-0.5)}-!{(-1.5,1.5)},
!{(0.5,-0.5)}-!{(1.5,0.5)},
!{(-0.5,0.5)}-!{(0.5,1.5)},
!{(0.5,-1.5)}-@(d,ur)!{(1.5,0.5)},
!{(0,0)}*{\scriptstyle{\bullet}},
!{(0.5,-1)}*{\scriptstyle{\circ}}
}.
\end{align*}
We show the computation of the two first terms of the left hand side,
$$
\partial
\xygraph{
!{<0pt,0pt>;<10pt,0pt>:<0pt,10pt>::},
!{(-1,1)}-!{(0,0)}-!{(1,1)},
!{(0,0)}*{-}
}
=
\xygraph{
!{<0pt,0pt>;<10pt,0pt>:<0pt,10pt>::},
!{(0,0)}-!{(0,-1)},
!{(0,0)}-!{(-1,1)},
!{(0,0)}-!{(1,1)},
!{(0,-1)}*{-},
}
+
\xygraph{
!{<0pt,0pt>;<10pt,0pt>:<0pt,10pt>::},
!{(0,0)}-!{(0,-1)},
!{(0,0)}-!{(-1,1)},
!{(0,0)}-!{(0,1)},
!{(0,0)}-!{(1,1)},
!{(0,-1)}-@(d,ur)!{(1,1)}
},
$$

\begin{align*}
\partial
\xygraph{
!{<0pt,0pt>;<10pt,0pt>:<0pt,10pt>::},
!{(-1,1)}-!{(0,0)}-!{(1,1)},
!{(0,0)}-!{(0,-1)},
!{(0,-0.5)}*{\scriptstyle{\bullet}},
!{(0,-1)}*{-}
}
=&
-
\xygraph{
!{<0pt,0pt>;<10pt,0pt>:<0pt,10pt>::},
!{(-1,1)}-!{(0,0)}-!{(1,1)},
!{(0,0)}-!{(0,-1)},
!{(0,-0.5)}*{\scriptstyle{\bullet}},
!{(0,-1)}*{-},
!{(-1,1)}*{\scriptstyle{+}}
}
-
\xygraph{
!{<0pt,0pt>;<10pt,0pt>:<0pt,10pt>::},
!{(-1,1)}-!{(0,0)}-!{(1,1)},
!{(0,0)}-!{(0,-1)},
!{(0,-0.5)}*{\scriptstyle{\bullet}},
!{(0,-1)}*{-},
!{(1,1)}*{\scriptstyle{+}}
}
=
\xygraph{
!{<0pt,0pt>;<10pt,0pt>:<0pt,10pt>::},
!{(-1,1)}-!{(0,0)}-!{(1,1)},
!{(0,0)}-!{(0,-1.5)},
!{(0,-1)}*{\scriptstyle{\bullet}},
!{(0,-1.5)}*{-},
!{(0,-0.5)}*{\scriptstyle{\times}}
}
=
-
\xygraph{
!{<0pt,0pt>;<10pt,0pt>:<0pt,10pt>::},
!{(-1,1)}-!{(0,0)}-!{(1,1)},
!{(0,0)}-!{(0,-1.5)},
!{(0,-0.5)}*{\scriptstyle{\bullet}},
!{(0,-1.5)}*{-},
!{(0,-1)}*{\scriptstyle{\times}}
}
+
\xygraph{
!{<0pt,0pt>;<10pt,0pt>:<0pt,10pt>::},
!{(-1,1)}-!{(0,0)}-!{(1,1)},
!{(0,0)}-!{(0,-1)},
!{(0,-0.5)}*{\scriptstyle{\circ}},
!{(0,-1)}*{-}
}
-
\xygraph{
!{<0pt,0pt>;<10pt,0pt>:<0pt,10pt>::},
!{(-1,1)}-!{(0,0)}-!{(1,1)},
!{(0,0)}-!{(0,-1)},
!{(0,-1)}*{-}
}
=
\xygraph{
!{<0pt,0pt>;<10pt,0pt>:<0pt,10pt>::},
!{(0.5,-0.5)}-!{(0.5,-1.5)},
!{(0.5,-0.5)}-!{(-1.5,1.5)},
!{(0.5,-0.5)}-!{(1.5,0.5)},
!{(-0.5,0.5)}-!{(0.5,1.5)},
!{(0.5,-1.5)}-@(d,ur)!{(1.5,0.5)},
!{(0,0)}*{\scriptstyle{\bullet}}
}
+
\xygraph{
!{<0pt,0pt>;<10pt,0pt>:<0pt,10pt>::},
!{(-1,1)}-!{(0,0)}-!{(1,1)},
!{(0,0)}-!{(0,-1)},
!{(0,-0.5)}*{\scriptstyle{\circ}},
!{(0,-1)}*{-}
}
-
\xygraph{
!{<0pt,0pt>;<10pt,0pt>:<0pt,10pt>::},
!{(-1,1)}-!{(0,0)}-!{(1,1)},
!{(0,0)}-!{(0,-1)},
!{(0,-1)}*{-}
}.
\end{align*}
Proceeding in the same manner with the remaining terms we obtain
\begin{align*}
\partial
\xygraph{
!{<0pt,0pt>;<10pt,0pt>:<0pt,10pt>::},
!{(0.5,-0.5)}-!{(0.5,-1.5)},
!{(0.5,-0.5)}-!{(-1.5,1.5)},
!{(0.5,-0.5)}-!{(1.5,0.5)},
!{(-0.5,0.5)}-!{(0.5,1.5)},
!{(0.5,-1.5)}-@(d,ur)!{(1.5,0.5)},
!{(0.5,-1)}*{\scriptstyle{\bullet}},
!{(0,0)}*{\scriptstyle{\bullet}}
}
=&
\xygraph{
!{<0pt,0pt>;<10pt,0pt>:<0pt,10pt>::},
!{(0.5,-0.5)}-!{(0.5,-1.5)},
!{(0.5,-0.5)}-!{(-1.5,1.5)},
!{(0.5,-0.5)}-!{(1.5,0.5)},
!{(-0.5,0.5)}-!{(0.5,1.5)},
!{(0.5,-1.5)}-@(d,ur)!{(1.5,0.5)},
!{(0.5,-1)}*{\scriptstyle{\circ}},
!{(0,0)}*{\scriptstyle{\bullet}}
}
-
\xygraph{
!{<0pt,0pt>;<10pt,0pt>:<0pt,10pt>::},
!{(0.5,-0.5)}-!{(0.5,-1.5)},
!{(0.5,-0.5)}-!{(-1.5,1.5)},
!{(0.5,-0.5)}-!{(1.5,0.5)},
!{(-0.5,0.5)}-!{(0.5,1.5)},
!{(0.5,-1.5)}-@(d,ur)!{(1.5,0.5)},
!{(0,0)}*{\scriptstyle{\bullet}}
}
+
\xygraph{
!{<0pt,0pt>;<10pt,0pt>:<0pt,10pt>::},
!{(0.5,-0.5)}-!{(0.5,-1.5)},
!{(0.5,-0.5)}-!{(-1.5,1.5)},
!{(0.5,-0.5)}-!{(1.5,0.5)},
!{(-0.5,0.5)}-!{(0.5,1.5)},
!{(0.5,-1.5)}-@(d,ur)!{(1.5,0.5)},
!{(0.5,-1)}*{\scriptstyle{\bullet}},
!{(0,0)}*{\scriptstyle{\circ}}
}
-
\xygraph{
!{<0pt,0pt>;<10pt,0pt>:<0pt,10pt>::},
!{(0.5,-0.5)}-!{(0.5,-1.5)},
!{(0.5,-0.5)}-!{(-1.5,1.5)},
!{(0.5,-0.5)}-!{(1.5,0.5)},
!{(-0.5,0.5)}-!{(0.5,1.5)},
!{(0.5,-1.5)}-@(d,ur)!{(1.5,0.5)},
!{(0.5,-1)}*{\scriptstyle{\bullet}}
}, \\
\partial
\xygraph{
!{<0pt,0pt>;<10pt,0pt>:<0pt,10pt>::},
!{(-0.5,-0.5)}-!{(-0.5,-1.5)},
!{(-0.5,-0.5)}-!{(-1.5,0.5)},
!{(-0.5,-0.5)}-!{(1.5,1.5)},
!{(0.5,0.5)}-!{(-0.5,1.5)},
!{(-0.5,-1.5)}-@(d,ur)!{(1.5,1.5)},
!{(-0.5,-1)}*{\scriptstyle{\bullet}},
!{(0,0)}*{\scriptstyle{\bullet}}
}
=&
\xygraph{
!{<0pt,0pt>;<10pt,0pt>:<0pt,10pt>::},
!{(-0.5,-0.5)}-!{(-0.5,-1.5)},
!{(-0.5,-0.5)}-!{(-1.5,0.5)},
!{(-0.5,-0.5)}-!{(1.5,1.5)},
!{(0.5,0.5)}-!{(-0.5,1.5)},
!{(-0.5,-1.5)}-@(d,ur)!{(1.5,1.5)},
!{(-0.5,-1)}*{\scriptstyle{\bullet}},
!{(0,0)}*{\scriptstyle{\circ}}
}
-
\xygraph{
!{<0pt,0pt>;<10pt,0pt>:<0pt,10pt>::},
!{(-0.5,-0.5)}-!{(-0.5,-1.5)},
!{(-0.5,-0.5)}-!{(-1.5,0.5)},
!{(-0.5,-0.5)}-!{(1.5,1.5)},
!{(0.5,0.5)}-!{(-0.5,1.5)},
!{(-0.5,-1.5)}-@(d,ur)!{(1.5,1.5)},
!{(-0.5,-1)}*{\scriptstyle{\bullet}}
}
+
\xygraph{
!{<0pt,0pt>;<10pt,0pt>:<0pt,10pt>::},
!{(-0.5,-0.5)}-!{(-0.5,-1.5)},
!{(-0.5,-0.5)}-!{(-1.5,0.5)},
!{(-0.5,-0.5)}-!{(1.5,1.5)},
!{(0.5,0.5)}-!{(-0.5,1.5)},
!{(-0.5,-1.5)}-@(d,ur)!{(1.5,1.5)},
!{(-0.5,-1)}*{\scriptstyle{\circ}},
!{(0,0)}*{\scriptstyle{\bullet}}
}
-
\xygraph{
!{<0pt,0pt>;<10pt,0pt>:<0pt,10pt>::},
!{(-0.5,-0.5)}-!{(-0.5,-1.5)},
!{(-0.5,-0.5)}-!{(-1.5,0.5)},
!{(-0.5,-0.5)}-!{(1.5,1.5)},
!{(0.5,0.5)}-!{(-0.5,1.5)},
!{(-0.5,-1.5)}-@(d,ur)!{(1.5,1.5)},
!{(0,0)}*{\scriptstyle{\bullet}}
}, \\
\partial
\xygraph{
!{<0pt,0pt>;<10pt,0pt>:<0pt,10pt>::},
!{(0,0)}-!{(0,-1)},
!{(0,0)}-!{(-1,1)},
!{(0,0)}-!{(0,1)},
!{(0,0)}-!{(1,1)},
!{(0,-1)}-@(d,ur)!{(1,1)},
!{(0,-0.5)}*{\scriptstyle{\bullet}}
}
=&
\xygraph{
!{<0pt,0pt>;<10pt,0pt>:<0pt,10pt>::},
!{(0,0)}-!{(0,-1)},
!{(0,0)}-!{(-1,1)},
!{(0,0)}-!{(0,1)},
!{(0,0)}-!{(1,1)},
!{(0,-1)}-@(d,ur)!{(1,1)},
!{(0,-0.5)}*{\scriptstyle{\circ}}
}
-
\xygraph{
!{<0pt,0pt>;<10pt,0pt>:<0pt,10pt>::},
!{(0,0)}-!{(0,-1)},
!{(0,0)}-!{(-1,1)},
!{(0,0)}-!{(0,1)},
!{(0,0)}-!{(1,1)},
!{(0,-1)}-@(d,ur)!{(1,1)}
}
+
\xygraph{
!{<0pt,0pt>;<10pt,0pt>:<0pt,10pt>::},
!{(-0.5,-0.5)}-!{(-0.5,-1.5)},
!{(-0.5,-0.5)}-!{(-1.5,0.5)},
!{(-0.5,-0.5)}-!{(1.5,1.5)},
!{(0.5,0.5)}-!{(-0.5,1.5)},
!{(-0.5,-1.5)}-@(d,ur)!{(1.5,1.5)},
!{(-0.5,-1)}*{\scriptstyle{\bullet}}
}
+
\xygraph{
!{<0pt,0pt>;<10pt,0pt>:<0pt,10pt>::},
!{(-0.5,-0.5)}-!{(-0.5,-1.5)},
!{(-0.5,-0.5)}-!{(-1.5,0.5)},
!{(-0.5,-0.5)}-!{(1.5,1.5)},
!{(0.5,0.5)}-!{(-0.5,1.5)},
!{(-0.5,-1.5)}-@(d,ur)!{(1.5,1.5)},
!{(0,0)}*{\scriptstyle{\bullet}}
}
+
\xygraph{
!{<0pt,0pt>;<10pt,0pt>:<0pt,10pt>::},
!{(0.5,-0.5)}-!{(0.5,-1.5)},
!{(0.5,-0.5)}-!{(-1.5,1.5)},
!{(0.5,-0.5)}-!{(1.5,0.5)},
!{(-0.5,0.5)}-!{(0.5,1.5)},
!{(0.5,-1.5)}-@(d,ur)!{(1.5,0.5)},
!{(0.5,-1)}*{\scriptstyle{\bullet}}
},
\end{align*}
which shows the wanted equality.
\eex
\section{Differential geometric interpretation}

Given a dg vector space $\fg$, it is well-known that structures of $L_\infty$-algebra on $\fg$ are in one-to-one correspondence to degree $1$ square-zero coderivations of the cofree graded cocommutative coalgebra without counit on $\fg[1]$. We consider its dual, the completed free graded commutative algebra without unit on $\fg[1]^*$, as the ideal of the basepoint in the ring of functions $\cO_{\fg[1]}=\prod_{n\geq 0}\odot^n \fg[1]^*$ on the formal graded pointed manifold $(\fg[1],0)$. We denote $\Der\cO_{\fg[1]}$ by $\cT_{\fg[1]}$ and refer to its elements as vector fields. A degree $1$ coderivation of $\odot^{\geq 1}\fg[1]$ corresponds, by dualizing, to a degree $1$ derivation of $\cO_{\fg[1]}$ which vanish at $0$. Hence, an $L_\infty$-structure on $\fg$ is given by a vector field $Q\in\cT^1\cO_{\fg[1]}$ such that $Q^2=0$ and $Q|_0=0$. Such vector fields are called homological.

Let now $(\fg,\{\ell_n\},\{q_n\})$ be a unimodular $L_\infty$-algebra. The collection $\{\tilde{q}_n\}$ assemble into a morphism $\odot^{\geq 1}\fg[1]\ra\k$. On the dual side, this determines a degree zero element $f$ in the ideal of $0$ in $\cO_{\fg[1]}$. We will want to express the equations satisfied by $\{\ell_n\}$ and $\{q_n\}$ in terms of $Q$ and $f$. This we will do using the divergence of the canonical volume form on $\fg[1]$ along $Q$. The volume forms are elements of the Berezinian of $\cT_{\fg[1]}$ which we now describe. For reference, see \cite{Bernstein1977a,Bernstein1977b,Leites1980,Manin1997}.

Let $A$ be a graded ring and $M$ a free $A$-module. If we choose a basis $\{e_\alpha\}$ of $M$ such that $|e_\alpha|$ is even for $1\leq \alpha\leq p$ and odd for $p+1\leq \alpha\leq n$, then the matrix of an even map $f:M\ra M$ is of the form
$$
\begin{pmatrix}
A & B \\
C & D
\end{pmatrix}
$$
where all entries in $A$ and in $D$ are even, while all entries in $B$ and in $C$ are odd. One can show that such a map is invertible if and only if $A$ and $D$ are both invertible, and moreover that, in this case, the Berezinian,
$$
\Ber(f)=\frac{\det(A-BD^{-1}C)}{\det D},
$$
is independent of the choice of basis. For any free $A$-module $M$, $\Ber M$ is a free rank one $A$-module. Every choice of basis $\{e_\alpha\}$ of $M$ determines a basis element $\Delta_x$ and in different bases of $M$ these are related by the Berezinian of the change of basis.

Choosing a basis of $\fg$, say $\{e_\alpha\}$, determines a basis $\{x^\alpha\}$ of $\fg[1]^*$, i.e.~a local coordinate system on the manifold $\fg[1]$. Here $x^\alpha$ is the basis element dual to $x_\alpha=s^{-1}e_\alpha$. This gives an isomorphism of algebras $\cO_{\fg[1]}\iso\k[[x^\alpha]]$ and of free $\cO_{\fg[1]}$-modules $\cT_{\fg[1]}\iso\oplus_\alpha\cO_{\fg[1]}\cdot\partial_\alpha$, where $\partial_\alpha=\partial/\partial x^\alpha$. Basis elements of $\Ber\cT_{\fg[1]}$ are called volume forms on $\fg[1]$. Hence, choosing local coordinates on $\fg[1]$ determines a constant volume form $\Delta_x$.

There is a well-defined action of $\cT_{\fg[1]}$ on $\Ber\cT_{\fg[1]}$, written $L_V(\rho)$ for $V\in\cT_{\fg[1]}$ and $\rho\in\Ber\cT_{\fg[1]}$, determined, in a set of coordinates $\{x^\alpha\}$, by the conditions
\begin{enumerate}
\item[i)] $L_V(g\rho)=V(g)\rho+(-1)^{|V||g|}g L_V(\rho)$,
\item[ii)] $L_{gV}(\rho)=(-1)^{|g||V|}L_V(g\rho)$ and
\item[iii)] $L_{\partial_\alpha}(\Delta_x)=0$ for all $\alpha$,
\end{enumerate}
see \cite{Leites1980}.

The \emph{divergence} of a vector field $V$ is defined by $L_V(\rho)=\div_\rho V\cdot\rho$. In coordinates $\{x^\alpha\}$, if $V=V^\alpha\partial_\alpha$ and $\rho=g\Delta_x$, then
\begin{align*}
L_V(\rho)=&L_{V^\alpha\partial_\alpha}(g\Delta_x)=(-1)^{|x^\alpha||V^\alpha|}L_{\partial_\alpha}(V^\alpha g\Delta_x)=(-1)^{|x^\alpha||V^\alpha|}\partial_\alpha(V^\alpha g)\Delta_x \\
=&\left((-1)^{|x^\alpha||V^\alpha|}\partial_\alpha V^\alpha g+V^\alpha\partial_\alpha g\right)\Delta_x=\left((-1)^{|x^\alpha||V^\alpha|}\partial_\alpha V^\alpha+V^\alpha\partial_\alpha g g^{-1}\right)\rho.
\end{align*}
Hence,
$$
\div_\rho V=(-1)^{|x^\alpha||V^\alpha|}\partial_\alpha V^\alpha+V^\alpha\partial_\alpha g g^{-1}.
$$

The fact that, for a unimodular $L_\infty$-algebra $(\fg,\{\ell_n\},\{q_n\})$,
\begin{align*}
\tr_{n+1}\tilde{\ell}_{n+1}=&\tr_{n+1}\tilde{\ell}^\alpha_{(\beta_1\ldots\beta_{n+1})}x_\alpha x^{\beta_1}\cdots x^{\beta_{n+1}} \\
=&(-1)^{|x^{\beta_{n+1}}|(|\tilde{\ell}_{n+1}|-|x^{\beta_{n+1}}|)}\tilde{\ell}^\alpha_{(\beta_1\ldots\beta_{n+1})}x^{\beta_{n+1}}x_\alpha x^{\beta_1}\cdots x^{\beta_n}
=\tilde{\ell}^\alpha_{(\beta_1\ldots\beta_n)\alpha}x^{\beta_1}\cdots x^{\beta_n} \\
\end{align*}
equals
\begin{align*}
-\sum_{i+j=n}\tilde{q}_{i+1}(1\id\otimes\cdots\otimes\tilde{\ell}_j)=-\sum\tilde{q}_{(\gamma_1\ldots\gamma_i|\alpha|}\tilde{\ell}^\alpha_{\beta_1\ldots\beta_j)}x^{\gamma_1}\cdots x^{\gamma_i}x^{\beta_1}\cdots x^{\beta_j}
\end{align*}
translates, for the associated vector field $Q$ and function $f$, into the equality
$$
\div_{\Delta_x}Q+Q(f)=0.
$$

Now, let $f$ be any degree zero function and $\rho=g\Delta_x$ a volume form. We have
\begin{align*}
\div_{e^f\rho} V=&(-1)^{|x^\alpha||V^\alpha|}\partial_\alpha V^\alpha+V^\alpha\partial_\alpha(e^f g)g^{-1}e^{-f}
=(-1)^{|x^\alpha||V^\alpha|}\partial_\alpha V^\alpha+V^\alpha\partial_\alpha f+V^\alpha\partial_\alpha g g^{-1} \\
=&\div_\rho V+V(f).
\end{align*}
From now on we let $\rho_f=e^f\Delta_x$, $\div=\div_{\Delta_x}$ and $\div_f=\div_{\rho_f}$. It follows from the above that $L_Q(\rho_f)=\div_f Q\cdot \rho_f=(\div Q+Q(f))\rho_f$. Thus, $(Q,f)$ is a unimodular $L_\infty$-structure if and only if $L_Q(\rho_f)=0$, i.e.~if and only if $\rho_f$ is $Q$-invariant.

\section{The deformation complex}

The graded vector space $\cT_{\fg[1]}$ is a graded Lie algebra with commutator bracket. Concerning $\div_\rho:\cT_{\fg[1]}\ra\cO_{\fg[1]}$ we have the following well-known lemma.

\ble\label{lem:sixone}
The equality
$$
\div_\rho[V_1,V_2]=V_1(\div_\rho V_2)-(-1)^{|V_1||V_2|}V_2(\div_\rho V_1)
$$
holds.
\ele

Given an $L_\infty$-algebra $(\fg,Q)$, $Q$ is a differential on $\cO_{\fg[1]}$ and $[Q,-]$ is a differential on $\cT_{\fg[1]}$.

\ble
If $(\fg,Q,f)$ is a unimodular $L_\infty$-algebra, then
$$
\div_f:\cT_{\fg[1]}\ra\cO_{\fg[1]}
$$
is a morphism of dg vector spaces.
\ele
\bpro
By Lemma \ref{lem:sixone} and the equality $\div_f Q=0$ we have
$$
\div_f([Q,V])=Q(\div_f V)-(-1)^{|V|}V(\div_f Q)=Q(\div_f V).
$$
\epro

\bde
Given a unimodular $L_\infty$-algebra $(\fg,Q,f)$, its \emph{deformation complex} is defined as
$$
C(\fg,Q,f)=\cone(\div_f)[-1],
$$
that is, $C(\fg,Q,f)=\cT_{\fg[1]}\oplus\cO_{\fg[1]}[-1]$ with differential
$$
\partial(V+sg)=[Q,V]-sQ(g)+s\div_f V.
$$
\ede

We can extend the Lie bracket on $\cT_{\fg[1]}$ to all of $C(\fg,Q,f)$ by letting $[V,sg]=(-1)^{|V|}sV(g)$ and $[sg_1,sg_2]=0$. With this convention, the differential of the deformation complex satisfies
\begin{align*}
\partial(V+sg)=&[Q,V]-sQ(g)+s\div_f V
=s\div_f V+[Q,V]+[Q,sg]
=s\div_f V+[Q,V+g],
\end{align*}
i.e.
$$
\partial=s\div_f+[Q,-].
$$

From Lemma \ref{lem:sixone} follows that $s\div_f$ is a degree $1$ derivation of the bracket on $C(\fg,Q,f)$. Since $[Q,-]$ is also a degree $1$ derivation, it follows that the deformation complex is a dg Lie algebra. The definition of the deformation complex is motivated by the following theorem.

\bth
Let $(\fg,Q,f)$ be a unimodular $L_\infty$-algebra. Pairs $(V,g)$ such that $(\fg,Q+V,f+g)$ is again a unimodular $L_\infty$-algebra are in one-to-one correspondence with Maurer-Cartan elements in $C(\fg,Q,f)$.
\eth
\bpro
A Maurer-Cartan element $V-sg$ of $C(\fg,Q,f)$ satisfies
\begin{align*}
\partial(V-sg)+&\frac{1}{2}[V-sg,V-sg]=s\div_f V+[Q,V]-[Q,sg]+V^2-[V,sg] \\
=&QV+VQ+V^2+s(\div_f V+Q(g)+V(g))=0,
\end{align*}
i.e.
\begin{align*}
QV+VQ+V^2=&0, \\
\div_f V+(Q+V)(g)=&0.
\end{align*}
On the other hand, $(Q+V,f+g)$ is a unimodular $L_\infty$-structure if and only if
\begin{align*}
(Q+V)^2=Q^2+QV+VQ+V^2=QV+VQ+V^2=&0, \\
\div_{f+g}(V+Q)=\div_f Q+\div_f V+(V+Q)(g)=\div_f V+(V+Q)(g)=&0.
\end{align*}
\epro

\section{The Characteristic Class}

We briefly address the question of when a given $L_\infty$-algebra $(\fg,Q)$ can be extended to a unimodular $L_\infty$-algebra. Lemma \ref{lem:sixone} implies that
\[
Q(\div Q)=\frac{1}{2}\div([Q,Q])=0,
\]
so that $\div Q$ determines a cohomology class $[\div Q]\in H^1(\cO_{\fg[1]},Q)$.
\bde
Given an $L_\infty$-algebra $(\fg,Q)$, the cohomology class $[\div Q]$ is called its \emph{characteristic class}.
\ede
\bpr
An $L_\infty$-algebra can be extended to a unimodular $L_\infty$-algebra if and only if its characteristic class vanishes.
\epr
\bpro
The characteristic class vanishes if and only if $\div Q=Q(-f)$ for some $f$. If so, then $\div_f Q=\div Q+Q(f)=0$ so that $(\fg,Q,f)$ is unimodular $L_\infty$.
\epro
\appendix

\section{Unimodular associative algebras and characteristic classes}

\subsection{Unimodular associative algebras}\hspace{1pt}

An (associative) algebra $(A,\mu)$ is said to be unimodular if, for all $x\in A$, \\
$\tr\mu(x,-)=\tr\mu(-,x)=0$. Unimodular algebras are algebras over the wheeled operad
\[
\UAss=\frac{\cF^\cal\left(
{\begin{xy}
(0,0);(-2,2)**@{-},
(0,0);(2,2)**@{-},
(0,0);(0,-2)**@{-},
(-2,3)*{\scriptstyle{1}},
(2,3)*{\scriptstyle{2}}
\end{xy}},
{\begin{xy}
(0,0);(-2,2)**@{-},
(0,0);(2,2)**@{-},
(0,0);(0,-2)**@{-},
(-2,3)*{\scriptstyle{2}},
(2,3)*{\scriptstyle{1}}
\end{xy}}
\right)}
{\left(
{\begin{xy}
(1,-1);(-3,3)**@{-},
(1,-1);(1,-3)**@{-},
(-1,1);(1,3)**@{-},
(1,-1);(3,1)**@{-}
\end{xy}}-
{\begin{xy}
(-1,-1);(3,3)**@{-},
(-1,-1);(-1,-3)**@{-},
(1,1);(-1,3)**@{-},
(-1,-1);(-3,1)**@{-}
\end{xy}},
{\begin{xy}
(0,0);(-2,2)**@{-},
(0,0);(2,2)**@{-},
(0,0);(0,-2)**@{-},
(0,-2);(-2,2)**\crv{(0,-3)&(-4,-2)&(-3,3)},
\end{xy}},
{\begin{xy}
(0,0);(-2,2)**@{-},
(0,0);(2,2)**@{-},
(0,0);(0,-2)**@{-},
(0,-2);(2,2)**\crv{(0,-3)&(4,-2)&(3,3)},
\end{xy}},
\right)},
\]
and since $\UAss^!=\Ass^\cal$ is wheeled Koszul \cite{Markl2007}, we can construct $\UAss_\infty$ explicitly.

\bpr
The minimal resolution of $\UAss$ is a free wheeled operad on generators
\[
\left\{
{\begin{xy}
(0,-1);(-4,3)**@{-},
(0,-1);(-2,3)**@{-},
(0,-1);(0,3)**@{-},
(0,-1);(2,3)**@{-},
(0,-1);(4,3)**@{-},
(0,-1);(0,-3)**@{-},
(0,4)*{\scriptstyle{n}}
\end{xy}}
\right\}_{n\geq 2}
\hspace{1cm}\text{and}\hspace{1cm}
\left\{
{\begin{xy}
(0,-2);(-6,2)**@{-},
(0,-2);(-4,2)**@{-},
(0,-2);(-2,2)**@{-},
(0,-2);(2,2)**@{-},
(0,-2);(4,2)**@{-},
(0,-2);(6,2)**@{-},
(0,-2)*{\scriptstyle{\blacktriangle}},
(-3,4)*{\scriptstyle{m}},
(3,4)*{\scriptstyle{n}}
\end{xy}}
\right\}_{m+n\geq 1}
\]
of degrees $2-n$ and $-m-n$ respectively. The second type of generator satisfies the following cyclic symmetry
\[
(-1)^{m+1}
{\begin{xy}
(0,-2);(-6,2)**@{-},
(0,-2);(-4,2)**@{-},
(0,-2);(-2,2)**@{-},
(0,-2);(2,2)**@{-},
(0,-2);(4,2)**@{-},
(0,-2);(6,2)**@{-},
(0,-2)*{\scriptstyle{\blacktriangle}},
\end{xy}}\zeta=
(-1)^{n+1}
{\begin{xy}
(0,-2);(-6,2)**@{-},
(0,-2);(-4,2)**@{-},
(0,-2);(-2,2)**@{-},
(0,-2);(2,2)**@{-},
(0,-2);(4,2)**@{-},
(0,-2);(6,2)**@{-},
(0,-2)*{\scriptstyle{\blacktriangle}},
\end{xy}}\xi=
{\begin{xy}
(0,-2);(-6,2)**@{-},
(0,-2);(-4,2)**@{-},
(0,-2);(-2,2)**@{-},
(0,-2);(2,2)**@{-},
(0,-2);(4,2)**@{-},
(0,-2);(6,2)**@{-},
(0,-2)*{\scriptstyle{\blacktriangle}},
\end{xy}},
\]
where $\zeta=(1\;2\;\ldots\;m)$ and $\xi=(m+1\;m+2\;\ldots\;m+n)$. The differential is defined by
\[
\partial
{\begin{xy}
(0,-1);(-4,3)**@{-},
(0,-1);(-2,3)**@{-},
(0,-1);(0,3)**@{-},
(0,-1);(2,3)**@{-},
(0,-1);(4,3)**@{-},
(0,-1);(0,-3)**@{-},
(0,4)*{\scriptstyle{n}}
\end{xy}}
=
\sum_{j=2}^{n-1}\sum_{i=0}^{n-j}(-1)^{ij+i+j+1+n}
{\begin{xy}
(-1,-5);(-1,1)**@{-},
(-1,-3);(-5,1)**@{-},
(-1,-3);(-3,1)**@{-},
(-1,-3);(-1,1)**@{-},
(-1,-3);(1,1)**@{-},
(-1,-3);(3,1)**@{-},
(1,1);(-3,5)**@{-},
(1,1);(-1,5)**@{-},
(1,1);(1,5)**@{-},
(1,1);(3,5)**@{-},
(1,1);(5,5)**@{-},
(-4,3)*{\scriptstyle{i}},
(1,7)*{\scriptstyle{j}},
\end{xy}}
\]
and
\begin{align*}
\partial
{\begin{xy}
(0,-2);(-6,2)**@{-},
(0,-2);(-4,2)**@{-},
(0,-2);(-2,2)**@{-},
(0,-2);(2,2)**@{-},
(0,-2);(4,2)**@{-},
(0,-2);(6,2)**@{-},
(0,-2)*{\scriptstyle{\blacktriangle}}
\end{xy}}=&
\sum_{i=2}^m\sum_{j=0}^{m-1}(-1)^{i+m+n+1}
{\begin{xy}
(2,-4);(-4,0)**@{-},
(2,-4);(-2,0)**@{-},
(2,-4);(0,0)**@{-},
(2,-4);(4,0)**@{-},
(2,-4);(6,0)**@{-},
(2,-4);(8,0)**@{-},
(-4,0);(-8,4)**@{-},
(-4,0);(-6,4)**@{-},
(-4,0);(-4,4)**@{-},
(-4,0);(-2,4)**@{-},
(-4,0);(-0,4)**@{-},
(2,-4)*{\scriptstyle{\blacktriangle}},
(-4,6)*{\scriptstyle{i}}
\end{xy}}
((-1)^{m+1}\zeta)^j
+\sum_{i=2}^n\sum_{j=0}^{n-1}
(-1)^{im+n}
{\begin{xy}
(0,-4);(-6,0)**@{-},
(0,-4);(-4,0)**@{-},
(0,-4);(-2,0)**@{-},
(0,-4);(2,0)**@{-},
(0,-4);(4,0)**@{-},
(0,-4);(6,0)**@{-},
(2,0);(-2,4)**@{-},
(2,0);(0,4)**@{-},
(2,0);(2,4)**@{-},
(2,0);(4,4)**@{-},
(2,0);(6,4)**@{-},
(0,-4)*{\scriptstyle{\blacktriangle}},
(2,6)*{\scriptstyle{i}}
\end{xy}}
((-1)^{n+1}\xi)^j \\
&+\sum_{i=0}^{m-1}\sum_{j=0}^{n-1}
{\begin{xy}
(0,-1);(-4,3)**@{-},
(0,-1);(-2,3)**@{-},
(0,-1);(0,3)**@{-},
(0,-1);(2,3)**@{-},
(0,-1);(4,3)**@{-},
(0,-1);(0,-3)**@{-},
(0,-3);(0,3)**\crv{(0,-5)&(2,-7)&(6,4)&(2,7)&(0,5)},
\end{xy}}
((-1)^{m+1}\zeta)^i((-1)^{n+1}\xi)^j.
\end{align*}
\epr

\subsection{Characteristic classes}\hspace{1pt}

An $A_\infty$-algebra, i.e.~a representation of $\Ass_\infty$ in a vector space $V$, corresponds to a degree one element
\[
\Gamma=\{\Gamma_n:V[1]^{\otimes n}\ra V[1]\}_{n\geq 1}
\]
in the weight graded Lie algebra
\[
C(V,V)=\{C^p_n(V,V)=\Hom^p(V[1]^{\otimes n},V[1])\}_{n\geq 1}
\]
equipped with the Gerstenhaber bracket, such that $[\Gamma,\Gamma]=0$. Consider next the right $C(V,V)$-module $\Cyc(V,\k)$, defined as
\[
\{\Cyc^p_{m,n}(V,\k)=\Hom^p((V[1]^{\otimes m})_{C_m}\otimes(V[1]^{\otimes n})_{C_n},\k)\}_{m+n\geq 1}.
\]
The module structure induces a Lie algebra structure on $\fg=C(V,V)\oplus\Cyc(V,\k)$ and the fact that $[\Gamma,\Gamma]=0$ implies that $\partial_\Gamma=[\Gamma,-]$ is a differential on $\fg$ such that $\Cyc(V,\k)$ is a subcomplex. Moreover, $\Gamma$ determines an element $\Gamma^\cal$ in $\Cyc^1(V,\k)$ as follows,
\[
\Gamma^\cal_{m,n}=\sum_{i=1}^m\sum_{j=1}^n\tr_{m+1}(\Gamma_{m+n+1})((-1)^{m+1}\zeta)^i((-1)^{n+1}\xi)^j.
\]
Here, $\zeta$ and $\xi$ are as above.

\bpr
For any $A_\infty$-algebra $(V,\Gamma)$ we have $\partial_\Gamma\Gamma^\cal=0$.
\epr
\bpro
For $m,n\geq 1$ this is Proposition-Definition $6.7.2$ in \cite{Markl2007}. For $m+n=1$ the proof is exactly the same, i.e.~expand the expressions $\partial\tr_1(\Gamma_n)$ and $\partial\tr_n(\Gamma_n)$ and skew symmetrize cyclically.
\epro
Hence $\Gamma$ determines a class $[\Gamma^\cal]$ in $H^1(\Cyc(V,\k),\partial_\Gamma)$ called the \emph{unimodular cyclic characteristic class}. The following proposition is proven in exactly the same way as Theorem $6.7.3$ in \cite{Markl2007}.
\bpr
A structure of $A_\infty$-algebra on $V$ can be extended to a structure of unimodular $A_\infty$-algebra if and only if the unimodular cyclic characteristic class vanishes.
\epr

\bibliographystyle{amsalpha}
\bibliography{umlie}

\end{document}